\documentclass[english,12pt,a4paper]{article}
\usepackage{adjustbox}
\usepackage{algorithm}
\usepackage{algpseudocode}
\usepackage{aligned-overset}
\usepackage{amsmath,amsthm}
\usepackage{authblk}
\usepackage[style=authoryear,maxbibnames=99,bibencoding=utf8,firstinits=true,backend=bibtex]{biblatex}
\usepackage{bm}
\usepackage{caption}
\usepackage{cases}
\usepackage[overload]{empheq}
\usepackage{float}
\usepackage[hmargin={26mm,26mm},vmargin={30mm,35mm}]{geometry}
\usepackage{graphicx}
\usepackage{hyperref}\usepackage{cleveref}
   \crefformat{section}{\S#2#1#3}
   \crefformat{subsection}{\S#2#1#3}
   \crefformat{subsubsection}{\S#2#1#3}
\usepackage{mwe}
\usepackage{newtxtext,newtxmath}
\usepackage{nicefrac}
\usepackage{subcaption}
\usepackage{tikz}\usetikzlibrary{shapes.geometric, arrows, positioning}
\usepackage{tikz-cd}
\usepackage{verbatim}
\usepackage{xcolor}
\usepackage{lipsum}
\usepackage{listings}
\usepackage{scalerel}

\newcommand{\email}[1]{\href{mailto:#1}{#1}}

\usepackage[normalem]{ulem}
\newcounter{corr}
\definecolor{violet}{rgb}{0.580,0.,0.827}
 \newcommand{\corr}[3]{\typeout{Warning : a correction remains in page \thepage}
   \stepcounter{corr}
 	            {\color{blue}\ifmmode\text{\,\sout{\ensuremath{#1}}\,}\else\sout{#1}\fi}
               {\color{red}#2}
               {\color{violet} #3}
 }

\newcommand{\partitle}[1]{
\\[3mm]
\underline{\textit{#1}}
\\[2mm]
}

\algnewcommand{\Initialize}[1]{%
  \State \textbf{Initialize:}
  \Statex \hspace*{\algorithmicindent}\parbox[t]{.8\linewidth}{\raggedright #1}
}

\theoremstyle{plain}
\newtheorem{theorem}{Theorem}

\theoremstyle{remark}
\newtheorem{remark}[theorem]{Remark}

\theoremstyle{definition}

\DeclareMathOperator*{\argmin}{arg\,min}

\newcommand{\Real}{\mathbb{R}}

\newcommand{\INT}{\mathrm{int}}
\newcommand{\EXT}{\mathrm{ext}}

\newcommand{\Th}{\mathcal{T}_h}
\newcommand{\Tho}{\mathcal{T}_{h,0}}

\newcommand{\Eh}{\mathcal{E}_h}

\newcommand{\ET}{\mathcal{E}_T}

\newcommand{\TE}{\mathcal{T}_E}

\newcommand{\Poly}[1]{\mathcal{P}^{#1}}

\newcommand{\Pone}[1]{\mathbb{P}^1(#1)}

\newcommand{\jumpg}[1]{[#1]_\gamma}

\newcommand{\norm}[2]{\|#2\|_{#1}}

\newcommand{\tnorm}[2]{|\kern-0.25ex|\kern-0.25ex|#2|\kern-0.25ex|\kern-0.25ex|_{#1}}

\newcommand{\E}{\mathcal{E}}
\newcommand{\T}{\mathcal{T}}
\newcommand{\V}{\mathcal{V}}

\newcommand{\Id}[1]{\mathcal{I}_{#1}}
\newcommand{\Ca}{\mathfrak{Ca}}

\bibliography{biblio}
\graphicspath{{./figures/}}

\title{Numerical simulation of biological capsules based on interface-fitted polygonal meshes}
\author{Daniele A. Di Pietro}
\author{Simon Mendez}
\author{Aurelio E. Spadotto}
\affil{IMAG, Univ. Montpellier, CNRS, Montpellier, France,
  \email{daniele.di-pietro@umontpellier.fr},
  \email{simon.mendez@umontpellier.fr},
  \email{aurelio-edoardo.spadotto@umontpellier.fr}
}

\begin{document}

\maketitle

\begin{abstract}
  In this work we introduce a numerical framework for the simulation of biological capsules based on meshes that adapt iteration-wise to the shape of the membrane.
  The discretization leverages a Hybrid High Order method (HHO) to solve a Stokes problem with a source term distributed on the interface between the internal and the external phases. 
  The HHO discretization supports meshes composed of general polygons, permitting to use fitted meshes that can be easily generated by cutting the elements intersected by the interface.
  Moreover, polygonal meshes allow arbitrarily accurate interface discretizations without introducing bulk degrees of freedom.
  The presentation of the simulation algorithm is completed with a proof-of-concept numerical test on a well-established benchmark: the deformation of an elastic capsule immersed in a shear flow. 
  \medskip\\
  \textbf{Key words:}
  Capsule Dynamics,
  Interface tracking,
  Interface Stokes problem,
  Fitted Methods,
  Polygonal Meshes,
  Hybrid High Order Method
\end{abstract}
\section{Introduction}
Capsules are physical systems that can be described as thin, soft membranes enclosing a fluid core.
Artificial capsules find applications in biotechnology, e.g. for drug design and the culture of artificial organs [\cite{barthes-bieselModelingMotionCapsules2011}];[\cite{lanzaCellEncapsulationTechnology1999a}]. 

Moreover, they stand as a mechanical model of living cells.
In numerical simulation of red blood cells, the latter are often treated as elastic capsules suspended in a fluid medium. 
The accurate description of capsules is still an open research front in biomechanics, and several models can be adopted as a starting point for a computational approach.
For a general overview on modeling aspects, the interested reader is referred to [\cite{vlahovskaVesiclesRedBlood2009}],[\cite{vlahovskaFlowDynamicsRed2013}] and, for a focus on the applications to red blood cells, to [\cite{pozrikidis2010computational}], [\cite{viallatDynamicsBloodCell2019}].

The flow of such fluidic particles is of particular interest for applications like microcirculation [\cite{Balogh:2017}], drug delivery [\cite{Liu:2019}], diagnostics [\cite{Krueger:2014,Pottier:2026}], or in fundamental understanding of particle or suspensions behavior, and simulation has emerged as a precious tool to complement microfluidic experiments [\cite{Iss:2019}].

One of the most common modeling strategies considers the membrane of the capsule as an elastic surface suspended in a flow which on the one hand transports it, and on the other hand receives a force originating from the elastic response of the membrane [\cite{breyiannisSimpleShearFlow2000, barthes1981time}].
Therefore, simulation of capsules is a particular  application of a broader class of problems of computational mechanics:
the numerical representation of fluid-structure interactions at immersed interfaces.

Although a rigorous classification of numerical methods for simulating capsule dynamics would go beyond the scope of this work,
a first broad distinction can be between approaches that discretize the volume domain and those that do not.
Indeed, artificial capsules or cells are often of a few microns in diameter, so that inertial effects
may often be neglected. The flow can then be described by the Stokes equations, which allow to express the
fluid-structure interaction problem in terms of boundary integrals, eliminating the need to mesh the bulk [\cite{pozrikidis1992boundary}]. The resulting equations may be solved numerically using the Boundary Element Method (BEM), which may be considered as the reference method for a small number of capsules in Stokes flows [\cite{Becic:2025,brebbia1978boundary,Dupont:2016,Peng:2011,Puthumana:2024}]. Apart from its limitation to linear equations, BEM also requires the surface mesh to be refined as the distance between the surface and other boundaries decreases and is not adapted to large suspension flows as it feature a dense matrix problem [\cite{pozrikidis1992boundary}].

As an alternative to BEM, a meshing of the domain is introduced to discretize flow equations.
Among mesh-based techniques, a distinction can be made between interface-fitted methods, in which the mesh is required to conform to the surface,
and unfitted methods, in which this constraint is removed.
In an interface-fitted approach, the interface is represented by a subset of faces of the background mesh, facilitating the discretization of transmission conditions [\cite{Kloppel:2011}].
As a drawback, the displacement of the interface over time typically requires a corresponding adaptation of the mesh.
The advantage of unfitted methods lies in having a fixed, possibly structured background mesh that does not need to follow the evolving surface, ths removing the need for mesh adaptation at each iteration.

One of the most popular unfitted methods is the Immersed Boundary Method (IBM) [\cite{peskinImmersedBoundaryMethod2002}].
In this framework, the action of the elastic forces on the fluid is represented via an effective volumetric source in the Navier Stokes equations.
Conversely, the velocity of the fluid is transmitted to the surface, whose position can be updated at the end of a time iteration.
As for the numerical solution of flow equations, several techniques can be considered.
Mainstream examples are the Finite-Difference and the Finite-Volume Methods [\cite{temam2024navier}].
An IBM algorithm coupled with a Finite-Volume solver for Navier-Stokes equations is used for the red blood cell simulations presented in [\cite{mendez2014unstructured,Balogh:2017, gai2026capturing}]. 
A main feature of the method is the fact that the solid surface is meshed.
Therefore, common discretizations of surface forces based on triangulations of the membrane ([\cite{charrier1989free}], [\cite{farutin20143d}]) are compatible with IBM.
As a drawback, the method requires the surface cell diameter to match that of the fluid domain, posing a constraint on the refinement on the mesh close to the surface.

For the sake of completeness, one may also mention other simulation techniques which are not based on solving the Stokes or the Navier-Stokes equations for the fluid flow, such as the lattice-Bolzmann-IBM technique [\cite{Krueger:2014,Liu:2019}], or particle methods such as the Dissipative Particle Dynamics [\cite{Fedosov:2014b}], to name a few. 

A different approach to the fluid/solid-surface interaction problem identifies the feedback of the membrane in terms of a transmission condition across the interface between the fluid phases.
In this sense, the simulation of capsules is directly linked to the numerical solution of problems featuring evolving interfaces.
The literature on the subject is large; a general overview can be found in [\cite[Sec. 16.10]{quarteroniNumericalApproximationPartial1994}].
Common techniques rely on an implicit representation of the surface by means of a scalar field defined on the fluid domain.
This approach avoids the complexity introduced by an explicit meshing of the surface.
Notable examples of unfitted low-order methods are the Volume Of Fluid method (VOF) [\cite{hirt1981volume}], the Level Set Method (LSM) [\cite{osher2001level}] and the Ghost Fluid Method (GFM) [\cite{fedkiwNonoscillatoryEulerianApproach1999}].
These techniques have been generalized to arbitrary order numerical schemes, in particular the Finite Element Method.
The family of unfitted Finite Element discretizations of interface problems includes several examples.
The XFEM method [\cite{belytschko1999arbitrary}] is based on enriching the trial/test space with discontinuous functions allowing the capture of jumps.
It has been adapted to interface Stokes equations in [\cite{chessa2003enriched}].
In the Cut-FEM method [\cite{burmanCutFEMDiscretizingGeometry2015}], jumps are represented thanks by duplicating the degrees of freedom over elements intersected by interface.
An application to fluid dynamics models is found in [\cite{massing2014stabilized}].
Still on Finite Elements, the $\phi$-FEM method [\cite{duprez2020phi}] ([\cite{duprez2023phi}] for applications to Stokes) allows a representation of jumps by integrating the level set function in the trial space.
The unfitted approach has also been generalized to hybrid schemes such as the Hybrid High Order Method [\cite{burmanUnfittedHybridHighorder2018}]; [\cite{burman2021unfitted}].

The objective of this work is to present a new approach to the numerical simulation of evolving interfaces.
The proposed framework leverages a polygonal discretization of the surface in the fashion of IBM and implements a fitted paradigm for the solution of an interface Stokes equation.
Interface-fitted meshes are simply generated by clipping a fixed background mesh along the discrete membrane.
The price to pay is the appearance of polygonal elements that exclude the adoption of conforming schemes.
This problem is solved by adopting the Hybrid High Order method (HHO [\cite{dipietroArbitraryorderCompactstencilDiscretization2014}];[\cite{dipietroHybridHighorderLockingfree2015}];[\cite{aghiliHybridizationMixedHighorder2015}]) for the solution of flow equations.
The HHO method is an arbitrary-order, nonconforming numerical scheme supporting general polygonal meshes.
Up to our knowledge, this is the first HHO solution of the interface Stokes problem and the first application of a polygonal scheme in the context of capsule simulation.

A feature of the solver proposed in this work is that the discretizaton of the membrane is substantially independent from the refinement level of the background mesh.
As a practical consequence, to reduce discretization errors linked to geometric approximation of the capsule, the surface can be refined while leaving to the background mesh intact.
This represents an appealing evolution with respect to mainstream approaches such as IBM and level-set-based methods, in which the accuracy on the membrane is limited by the discretization of the background.

The key passage in modeling capsule suspensions involves predicting the response of a single capsule to viscous forces exerted by the surrounding fluid.
In real-world contexts, elastic membranes are typically described by nonlinear constitutive laws, and capsules often experience very large deformations in typical flow configurations.
This makes capsule simulation a challenging problem where numerical simulation is a valuable tool for validating membrane models and exploring flow conditions not accessible to analytical treatment.
For validating numerical solvers, one particular flow configuration has emerged as a standard benchmark: the evolution of a spherical capsule freely suspended in a linear shear flow, to demonstrate the ability of the solver to accurately predict the maximal deformation of a capsule as a function of its membrane relative stiffness.
Therefore, this is the target test for a numerical simulator.

The rest of the article is structured as follows.
In Section~\ref{sec:pm-model}, the model of capsule dynamics is recalled.
In Section~\ref{sec:pm-simulator}, the design of the solver is exposed in detail. 
Section~\ref{sec:num-tests} is dedicated to numerical tests.
In particular, Section~\ref{subsec:test.shear} presents the predictions obtained for the capsule in shear flow benchmark.
The simulation output is successfully compared with previously published results obtained for the same test case.
Implementation aspects are discussed in Appendix~\ref{app.sec:polymemb}.
\section{Two-dimensional Capsule Model}\label{sec:pm-model}
\begin{figure}
	\centering
	\includegraphics[width=0.5\textwidth]{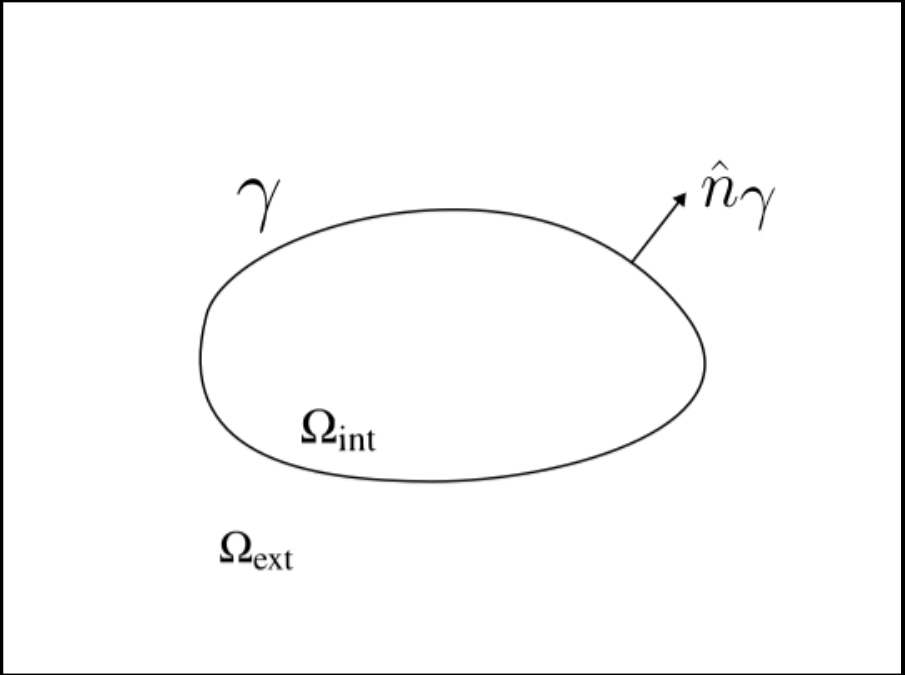}
    \caption{Representation of the 2D capsule described in Section~\ref{sec:pm-model}; the membrane $\gamma$ is described as a fluid/fluid interface with no thickness.}
    \label{fig:sketch_simple}
\end{figure}

Let $\Omega$ be an open bounded, polytopal domain of $\Real^2$.
The membrane of the capsule is represented by a closed, time-dependent curve $\gamma(t)$.
At any time, $\gamma (t)$ partitions $\Omega$ into the pair of open subdomains $\Omega_\INT$, $\Omega_\EXT$ (time dependence omitted for readability), representing respectively the internal and external fluid phases (see Fig.~\ref{fig:sketch_simple}).
Given a closed curve $\gamma\subset\Real^2$, let $\eta$ be the arc-length parametrization orienting $\gamma$ positively, such that the position of a point on $\gamma$ is provided by the map $\gamma(\eta)$ (to avoid the proliferation of symbols, we identify the support with the parametrization map using the same symbol $\gamma$).
The unitary tangent vector $\hat{t}_\gamma$, the unitary outward normal vector $\hat{n}_\gamma$ and the curvature $k$ are expressed by the relationships
\[
\hat{t}_\gamma = \frac{d \gamma}{d \eta}, \qquad -\hat{k}_\gamma = - k \hat{n}_\gamma = \frac{d \hat{t}_\gamma}{d \eta}.
\]
With $m\geq 0$, let $\mathbb{T}^m$ be the space of $m$-ranked tensors on $\Real^2$ ($\mathbb{T}^0 = \Real$, $\mathbb{T}^1 = \Real^2$).
Given a smooth field $z$ on $\gamma$ with values in $\mathbb{T}^m$, define the intrinsic gradient as the field of $(m+1)$-ranked tensors
\[
\nabla_\gamma z = \frac{d z}{d\eta} \otimes \hat{t}_\gamma,
\]
where the tensor product $(\cdot \otimes \cdot): \mathbb{T}^m\times \mathbb{T}^n \to \mathbb{T}^{m+n}$ is such that $ (a\otimes b)_{i_1,\dots, i_m,j_1,\dots,j_n} = a_{i_1,\dots,i_m} b_{j_1,\dots,j_n}$.
Consider also the intrinsic divergence and Laplacian operators
\[
\textrm{div}_\gamma z = \textrm{tr} (\nabla_\gamma z) = \frac{dz }{d\eta} \cdot \hat{t}_\gamma,
\qquad
\Delta_\gamma z = \textrm{div}_\gamma (\nabla_\gamma z) = \frac{d^2 z}{d \eta^2}.
\]
The membrane is assumed subjected to a force depending exclusively on its geometry and the state of deformation with respect to an undeformed configuration (this would not be the case for viscoelastic membrane for instance).
Force can be represented in terms of the curvilinear force density~$f_\gamma(\eta)$.
In this work, we consider two origins of force: strain and bending, which covers the applications of 2D capsules, vesicles and red blood cells [\cite{vlahovskaVesiclesRedBlood2009}].
To describe the elastic behaviour of the membrane, introduce a reference undeformed configuration $\gamma_0$ with arc-length $\eta_0$.
The current relative extension of the curve element is represented by the strain function
\[
\epsilon (\eta) = \frac{d \eta}{d \eta_0} - 1.
\]
Call $T(\eta)$ the scalar tension along $\gamma$. In this work we consider the linear elastic law
\[
T (\eta)= \lambda_{\rm{str}} \epsilon(\eta),
\]
where $\lambda_{\rm{str}}$ is the Young modulus of the membrane.
The resulting elastic force density reads
\begin{equation}\label{eq:f.str}
f_{\rm{str}} = \rm{div}_{\gamma} \big( T \, \hat{t}_\gamma \otimes \hat{t}_\gamma \big).
\end{equation}
The second source of force is linked by the principle of virtual works to the bending energy:
\begin{equation}\label{eq:bnd.energy}
E (\gamma)= \frac{\lambda_{\rm{bnd}}}{2} \int_\gamma k^2 d \eta,
\end{equation}
where $\lambda_{\rm{bnd}}$ is the bending modulus of the membrane.
The bending force density $f_{\rm{bnd}}(\eta)$ is defined such that, for any virtual displacement $\delta \gamma$, it holds
\[
\int_\gamma (f_{\rm{bnd}} \cdot \delta \gamma) d \eta = \delta E  = \frac{d}{d \varepsilon } \bigg( E\big(\gamma + \varepsilon \delta \gamma \big) \bigg)_{\varepsilon=0}.
\]
The above relationship leads to the following explicit expression:
\begin{equation}\label{eq:f.bnd}
   f_{\rm{bnd}} = \lambda_{\rm{bnd}} \bigg( \frac12 k^3  + \Delta_\gamma k \bigg) \hat{n}_\gamma.
\end{equation}
Putting together~\eqref{eq:f.str} and~\eqref{eq:f.bnd}, the total force density acting on $\gamma$ reads:
\[
    f_\gamma = f_{\rm{bnd}} + f_{\rm{str}}.
\]

\noindent
The membrane is advected by a velocity field $v_\gamma(\eta)$, such that $\gamma$ evolves according to the relationship $ \dot \gamma (\eta) = v_\gamma (\eta)$.
Velocity is imposed to the membrane by the surrounding fluid, which occupies the volume $\Omega_{\INT} \cup \Omega_{\EXT}$.
The fluid phase is described by a pressure $~{p:\Omega_\INT\cup\Omega_\EXT\to\Real}$ and velocity $v_{\Omega}:\Omega_\INT\cup\Omega_\EXT\to\Real^2$.
The flow is assumed incompressible, namely
\[
\nabla \cdot v_\Omega = 0 \qquad \text{in}\, \Omega.
\]
As a consequence, using $\jumpg{\phi}\coloneq \phi|_{\partial\Omega_{\INT}} - \phi|_{\partial\Omega_{\EXT}}$ to refer to the jump of a field $\phi$ across $\gamma$,
\[
\jumpg{v_\Omega} \cdot \hat{n}_\gamma = 0.
\]
Additionally, we consider the no-slip condition
\[
\jumpg{v_\Omega }\cdot \hat{t}_\gamma = 0.
\]
Since $v_\Omega$ is continuous in each component across $\gamma$, $v_\Omega$ can be extended to $\Omega$ and $v_\gamma$ can be taken as the restriction $(v_\Omega)_\gamma$.
With this premise, from now on we will write only $v$, without making the distinction between fluid velocity and membrane velocity.
The hydrodynamical force acting on a portion of fluid is assigned in terms of the tensor field $\sigma:\Omega_\INT\cup\Omega_\EXT\to \Real^{2\times2}$.
For a Newtonian fluid, the latter is described as a function of the couple $(p,v)$ as
\begin{equation}\label{eq:stress}
	\sigma = 2\mu \nabla_s v- p \Id{} = \mu(\nabla v + \nabla^T v) - p \Id{},
\end{equation}
where $\mu:\Omega_\INT\cup\Omega_\EXT\to\Real^2$ represents the dynamic viscosity of the fluid, possibly discontinuous across the interface,
$\nabla_s$ is the symmetric part of the gradient
and $\Id{}$ is the identity.
The velocity of the fluid is determined introducing a force balance.
At any time, for any volume region $X\subset\Omega$ with outward normal vector $\hat{n}_X$, the resulting force acting on the fluid occupying $X$ is zero.
If $X$ contains a section of $\gamma$, the internal force acting on the membrane section is transferred to the fluid.
This leads to the condition
\[
F(X) = \int_{X} \nabla \cdot \sigma + \int_{X \cap \gamma} f_\gamma = 0, \qquad \forall X\subset\Omega.
\] 
Since $X$ is arbitrary, the above integral condition is equivalent to the pair of local conditions
\[
   \left\{
   \begin{aligned}
      \nabla \cdot \sigma &= 0 &  \qquad \text{in} \,\Omega_\bullet,\, \bullet\in\{\INT,\EXT\},
      \\
      \jumpg{\sigma}\cdot \hat{n}_\gamma + f_\gamma &= 0 & \qquad\text{on}\,\gamma.
   \end{aligned}
   \right.
\]
To ensure uniqueness, given $g:\partial\Omega\to\Real^2$, we require
\[v = g \quad \text{on}\,\partial\Omega, \qquad \int_{\Omega} p = 0.\]
In summary, the membrane dynamics problem can be formulated as follows:
given an initial closed curve $\gamma_0\subset\Omega$, determine, for $t>0$, the curve $\gamma(t)$ such that
\begin{equation}
    \left\{
    \begin{aligned}
        \gamma(0) &= \gamma_0 &\\
        \dot \gamma & = v_\gamma &  \qquad &\forall t >0,
    \end{aligned}
    \right.
\end{equation}
where $v_\gamma = v|_\gamma$ and, for any $t$, $v:\Omega\to\Real^2$, $p:\Omega\to\Real$ are solution to the interface Stokes problem
\begin{equation}\label{eq:stokes-pb}
   \left\{
      \begin{aligned}
      \nabla \cdot \left( \mu(\nabla v + \nabla^Tv)-p \Id{}\right) &= 0 &  \qquad \text{in} \,&\Omega_{\bullet},\,\bullet\in\{\INT,\EXT\},
      \\
      \nabla \cdot v &= 0  & \qquad\text{in} \, &\Omega_{\bullet}, \,\bullet\in\{\INT,\EXT\},
      \\
      \jumpg{\mu(\nabla v + \nabla^Tv)-p \Id{} }\cdot \hat{n}_\gamma +f_\gamma &=0 & \qquad\text{on}\,&\gamma(t),
      \\
      \int_\Omega p& = 0, & &
      \\
      v &= g & \qquad \text{on} \,&\partial\Omega.
   \end{aligned}\right.
\end{equation}
We stress the fact that in problem~\eqref{eq:stokes-pb} there is no explicit dependence on time, and the solution $(p,v)$ only depends on the current configuration $\gamma(t)$.
\section{Simulation of Capsule Flow}\label{sec:pm-simulator}
The following Section is structured as follows:
in Section~\ref{subsec:intf.disc} the discretization of the capsule surface is discussed in detail;
in Section~\ref{subsec:hho-stokes} an HHO scheme for the interface Stokes problem is presented and the main convergence results are recalled with suitable references;
in Section~\ref{subsec:simul.algorithm} the time iteration scheme is presented together with a discussion on a volume correction strategy.

\subsection{Membrane discretization}\label{subsec:intf.disc}
The goal of this section is to illustrate the discrete model of the membrane at the base of the numerical simulation and a method to compute a discrete membrane force density.  
We wish to represent the membrane as a time-dependent polygon whose vertices have a time-dependent position.
If $\Gamma$ is a polygon, call $\mathcal{E}_{\Gamma}$ the set collecting the edges of $\Gamma$, $\V_\Gamma$ the set collecting its vertices and, for $E\in\E_\Gamma$, let $\Poly{k}(E, \mathbb{T}^m)$ be the space of $\mathbb{T}^m$-valued polynomials up to degree $k$ over $E$.
Define the space of broken polynomials
\[
\Poly{k}(\Gamma,\mathbb{T}^m) \coloneq \left\{ r: \Gamma\to \mathbb{T}^m \quad:\quad r|_E  \in \Poly{k}(E,\mathbb{T}^m) \quad \forall E\in\E_\Gamma \right\}.
\]
Consider the space $\mathbb{P}^k(\Gamma,\mathbb{T}^m)\coloneq C^0(\Gamma,\mathbb{T}^m) \cap \Poly{k}(\Gamma,\mathbb{T}^m)$,
where $C^0(\Gamma,\mathbb{T}^m)$ is the space of $\mathbb{T}^m$-valued continuous functions on $\Gamma$.
In this work we will make use of the spaces $\Pone{\Gamma, \mathbb{T}^m}$ and $\Poly{0}(\Gamma,\mathbb{T}^m)$, characterized respectively by degrees of freedom associated to $\V_\Gamma$ and $\E_\Gamma$.
\\[3mm]
We start by considering a closed polygonal chain $\Gamma_0$ approximating $\gamma_0$.
The current configuration of the membrane is given by a time-dependent transformation $\Gamma(t)$.
The transformation is supposed to be continuous and, at any time, the current configuration should also be a polygon.
In such a way it is possible to obtain the displacement of the discrete membrane by moving its vertices and reconnecting them by straight lines.
In other words, $\forall t>0$, $\Gamma(t)\in \mathbb{P}^1(\Gamma_0, \Real^2)$.
Notice then that $\hat{t}_\Gamma \in \mathcal{P}^0(\Gamma, \Real^2)$.
We look for a discrete membrane force density in $\mathbb{P}^1(\Gamma, \Real^2)$.
Given a closed curve $\gamma=\gamma(\eta)$ and a couple of tensor fields $p,q: \gamma \to \mathbb{T}^m$, consider the product
\[
(p,q)_\gamma \coloneq \int_\gamma p(x) \cdot q(x) d\eta.
\]
In order to introduce a discrete version of the elastic law~\eqref{eq:f.str}, we start by considering the integration by parts formula
\begin{equation}\label{eq:div-ipp}
(\textrm{div}_\gamma \sigma, z)_\gamma = -(\sigma, \nabla_\gamma z)_\gamma - (\sigma, z\otimes k \hat{n}_\gamma)_\gamma,
\end{equation}
valid when $\sigma$ and $z$ are sufficiently regular tensor field on $\gamma$, with respective rank $(m+1)$ and an $m$.
If $\sigma = (\sigma \cdot \hat{t}_\gamma)\otimes \hat{t}_\gamma$, the second term in~\eqref{eq:div-ipp} vanishes.
For example, when $\sigma = f_{\textrm{str}}$ or $\sigma = \Delta_\gamma u$ for some tensor field $u$, we get, respectively,
\begin{subequations}
\begin{align}
(f_{\textrm{str}}, z)_\gamma &= -(T \,\hat{t}_\gamma \otimes \hat{t}_\gamma, \nabla_\gamma z)_\gamma, \label{eq:ipp-f-str}
\\
(\Delta_\gamma u , z)_\gamma &= -( \nabla_\gamma u, \nabla_\gamma z)_\gamma. \label{eq:ipp-curv}
\end{align}
\end{subequations}
Projecting\eqref{eq:ipp-f-str} onto $\Pone{\Gamma, \Real^2}$, we introduce the discrete elastic force density $~{f_{\textrm{str},h} \in \mathbb{P}^1(\Gamma,\Real^2)}$ such that
\begin{equation}\label{eq:f.str.h}
(f_{\textrm{str,h}}, z)_\Gamma = - ( T\, \hat{t}_\Gamma\otimes\hat{t}_\Gamma, \nabla_\Gamma z)_\Gamma \qquad \forall z\in \mathbb{P}^1(\Gamma,\Real^2).
\end{equation}
\begin{remark}[Comparison with the literature]
The present approach matches the spring model adopted among others by [\cite{mendez2014unstructured}].
Consider a virtual displacement $\delta V\in\Real^2$ of a vertex $V$.
The latter corresponds to a the perturbation $\delta x = \varphi_V \delta V \in \Pone{\Gamma,\Real^2}$, where $\varphi_V\in\Pone{\Gamma,\Real}$ is the nodal basis hat function of vertex $V$.
The corresponding nodal force $F_V$ applied to $V$ is defined as
\[
F_V \cdot \delta V = (f_{\textrm{str},h}, \delta x)_\Gamma,
\]
where $(f_{\textrm{str},h}, \delta x)_\Gamma$ is the virtual work corresponding to $\delta x$.
Using $(-)$ and $(+)$ to label entities associated to the 2 the edges touching $V$ and applying~\eqref{eq:ipp-f-str}, we get
\[
(f_{\textrm{str},h}, \delta x )_\Gamma = - (T\, \hat{t}_\Gamma\otimes\hat{t}_\Gamma, \nabla_\Gamma \delta x)_\Gamma = (T_{(+)} \hat{t}_{(+)} - T_{(-)} \hat{t}_{(-)}) \cdot \delta V,
\]
whence
\[
F_V = T_{(+)} \hat{t}_{(+)} - T_{(-)} \hat{t}_{(-)}.
\]
Therefore, $f_{\textrm{str},h}$ expresses the same information of a system of forces applied at nodes in which each edge applies the tension $T_E$ on the neighbouring nodes.
\end{remark}
Projecting~\eqref{eq:ipp-curv} onto $\Pone{\Gamma, \Real^2}$, we introduce the discrete intrinsic Laplace operator $\Delta_{\Gamma}$ such that, for $~{ u\in\Pone{\Gamma, \Real^2}}$, $~{\Delta_\Gamma u \in \Pone{\Gamma, \Real^2}}$ satisfies
\begin{equation}\label{eq:disc-delta}
(\Delta_\Gamma u, z)_\Gamma = - ( \nabla_\Gamma u , \nabla_\Gamma z)_\Gamma \qquad \forall z\in \mathbb{P}^1(\Gamma,\Real^2).
\end{equation}
The bending force density is discretized as $f_{\textrm{bnd},h}\in \Pone{\Gamma,\Real^2}$ such that
\begin{equation}\label{eq:f.bnd.h}
(f_{\rm{bnd},h})_V \coloneq \lambda_{\rm{bnd}} \bigg( \|(\Delta_{\Gamma} \Gamma)_V\|^3  + \frac12 (\Delta_{\Gamma}^2\Gamma)_V \bigg) \qquad \forall V\in\V_\Gamma.
\end{equation}

\begin{remark}[Comparison with the literature]
The discrete intrinsic Laplace operator can be applied to obtain a discrete curvature in $\Pone{\Gamma,\Real^2}$. 
Let $V\in\V$, and use $(-)$, $(+)$ to label entities adjacent to $V$ from both sides.
Letting $\hat{k}_\Gamma = \Delta_{\Gamma} \Gamma$, from~\eqref{eq:disc-delta} we get, for all $ V \in \V$,
\[
\frac{1}{6}|E^{(+)}|\hat{k}_{(-)} + \frac{1}{3}(|E^{(-)}| + |E^{(+)}|)\hat{k}_V + \frac{1}{6}|E^{(+)}|\hat{k}_{(+)} = - |E^{(-)}|^{-1}(V- V^{(-)}) + |E^{(+)}|^{-1}(V^{(+)} - V),
\]
which, applying mass-lumping, reduces to
\[
\frac12 (|E^{(-)}|+|E^{(+)}|)  \hat{k}_V = \hat{t}_{(+)} - \hat{t}_{(+)},
\]
retrieving a well-established finite difference approximation of vector curvature.
\end{remark}
To conclude, summing together~\eqref{eq:f.str.h} and~\eqref{eq:f.bnd.h}, the total force density supported by the discrete membrane is
\begin{equation}
f_{\Gamma} = f_{\textrm{str},h} + f_{\textrm{bnd},h} \in \Pone{\Gamma,\Real^2}.
\end{equation}
\subsection{HHO approximation of flow velocity}\label{subsec:hho-stokes}
\begin{figure}[H]
	\centering
	\includegraphics[width=0.5\textwidth]{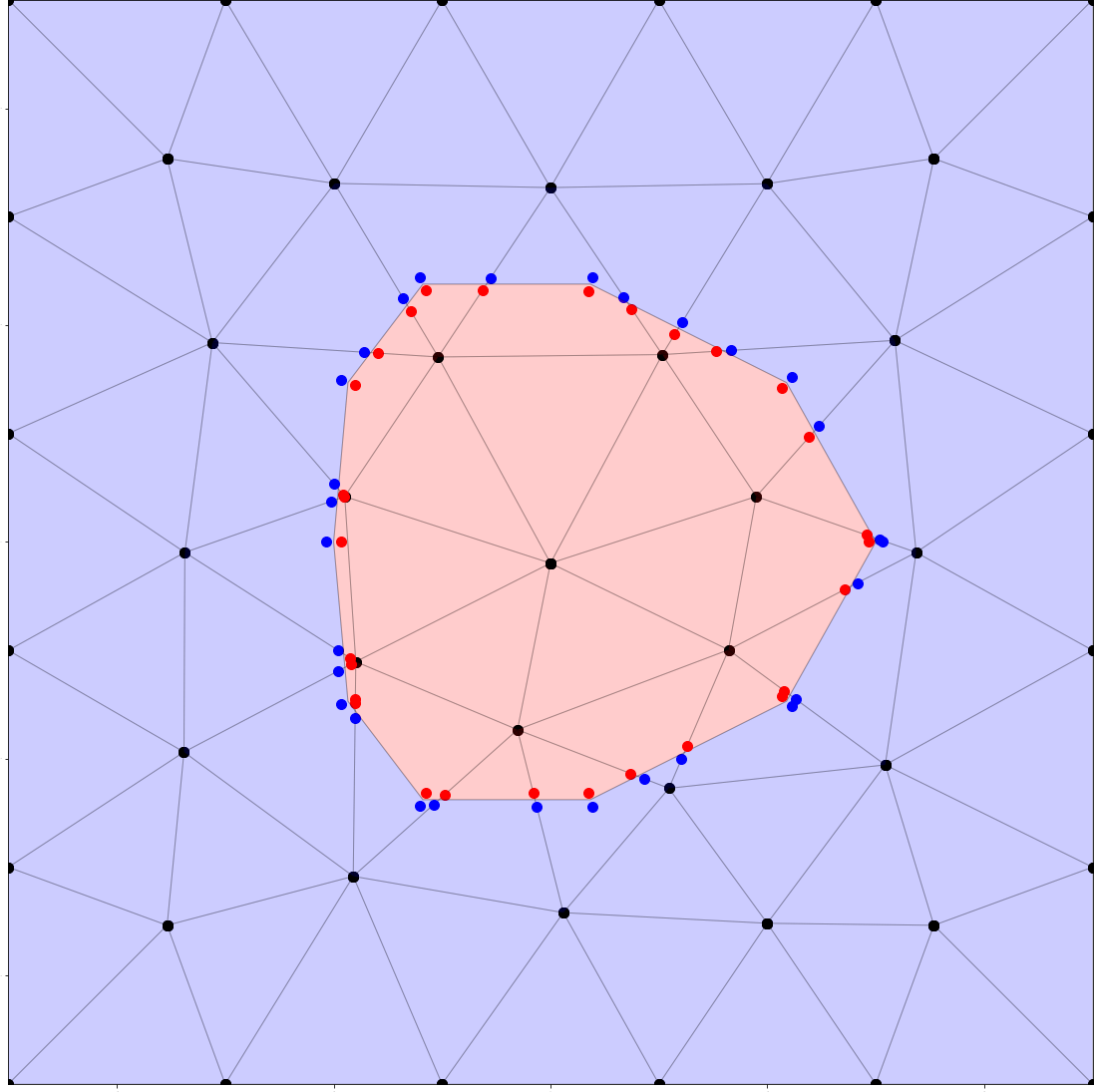}
    \caption{An example of admissible mesh compatible with an internal interface $\Gamma$.
    The discrete membrane $\Gamma$ is highlighted with colored double dots at the position of the nodes of $\V_\Gamma$.
    Elements are either internal (red) or external (blue), with no intersection with the interface but at the boundary.
    A compatible mesh can be generated starting from a standard triangular mesh by clipping the elements intersected by $\Gamma$.
    The resulting mesh features generic polygons spawning from cut triangles, motivating the introduction of a polygonal method.}
    \label{fig:compatible.mesh}
\end{figure}

At a given time $t>0$, the capsule boundary is discretized as a polygonal chain $\Gamma(t)$, supporting a force density $f_{\Gamma}\in\Pone{\Gamma,\Real^2}$ computed as described in Section~\ref{subsec:intf.disc}.
The velocity advecting the capsule is determined as the solution of problem~\eqref{eq:stokes-pb}.
As a constraint to design a numerical scheme, we require the domain $\Omega$ to be discretized by a mesh compatible with the membrane $\Gamma$.
This means essentially that all elements fall entirely only on one side of the membrane 
and that edges of the membrane are found at the interface between two elements (see Fig.~\ref{fig:compatible.mesh}).
In this sense, we wish to adpot a fitted approach, which provides a natural way to account for
the presence of an internal interface,
since the connectivity between interface edges and bulk elements is explicitely encoded in the mesh.
With moving membranes, such a requirement is demanding
since a new mesh must be generated after every displacement of $\Gamma$.
This can become particularly expensive for methods such as the FEM,
that has strict requirements on the shape of elements.
With the goal of preserving FEM-like arbitrary order accuracy whilst avoiding heavy mesh manipulation,
we opt for a high-order numerical scheme supporting generic polygonal meshes,
the Hybrid High Order Method.
This flexibility allows the mesh to conform to the moving interface through a simple clipping procedure applied to a fixed background mesh
An implementation of such a procedure procedure is illustrated in Appendix~\ref{app.sec:polymemb}.
Importantly, clipping generates polygonal elements close to the interface.

The method seeks an approximate velocity in a space featuring polynomial unknowns associated with both elements and edges.
The numerical scheme is based on two key ingredients: local reconstruction operators enabling $(k+1)$-order consistency errors derived from $k$-degree unknowns
and a stabilisation mechanism based on penalizing high-order differences.
An exhaustive theoretical study of the scheme for problems without membranes can be found in [\cite[Chapter 8]{dipietroHybridHighOrderMethod2020}].
Therefore, in the rest of the present section, we focus on the statement of the discrete problem, recalling briefly the main theoretical results and highlighting differences when necessary.
\noindent
\underline{\textit{Discrete Setting}}
\\[2mm]
Assume the domain $\Omega$ to be a polygon and introduce a mesh $(\T_h, \E_h)$, where the set of elements $\T_h$ is a collection of mutually disjoint polygons forming a partition of $\Omega$
and $\E_h$ is the set of element edges, partitioned into internal and boundary edges.
For an element $T\in\Th$, let $\E_T$ be the set of edges on its boundary.
Similarly, given $E\in\E_h$, let $\TE$ collect the two elements sharing $E$.
For $E\in\E_h$, let $n_E$ be a unit vector normal to $E$ chosen with arbitrary orientation;
for any couple $T,E$ with $E\in\E_T$, define the orientation $\omega_{TE}\in\{-1, 1\}$ such that $\omega_{TE}n_E$ is pointing out of $T$.
We assume the mesh to be compatible with $\Gamma$, that is, $\E_\Gamma\in\E_h$.
With $Y\in\Th \cup\Eh$, let $\Poly{k}(Y)$ be the space of scalar polynomials defined on $Y$
and let $\Poly{k}(Y,\Real^2)$ be the space of vector-valued polynomials.
\partitle{Discrete Spaces}
We introduce two spaces of polynomial degrees of freedom to approximate respectively velocity and pressure.
Let $k\geq1$ be the polynomial degree of the scheme.
For the velocity unknown, we introduce the space
\[
\underline{U}^k_h = \Big\{ \underline{w}_h = ((w_T)_{T\in\Th}, (w_E)_{E\in \Eh}): w_T\in \Poly{k}(T,\Real^2)\, \forall T\in\T_h,  w_E \in \Poly{k}(E,\Real^2) \forall E\in\E_h \Big\},\\
\]
featuring polynomial unknowns on both elements and edges.
With $T\in \mathcal T_h$, let $\underline{U}^k_T$ be the restriction of $\underline{U}^k_h$ to the entities attached to $T$.
Given $\underline{w}_h\in\underline{U}^k_h$, let $\underline{w}_T \coloneq (w_T, (w_E)_{E\in\ET})$.
The enrichment of the space with edge degrees of freedom has two consequences.
On the one hand, it allows to increase the order of approximation inside the elements by introducing local renconstructions.
For the Stokes problem, the velocity reconstruction operator ${r^{k+1}_T: \underline{U}^k_T \to \Poly{k+1}(T,\Real^2)}$
is defined such that
\begin{equation}
\begin{aligned}
	(\nabla r^{k+1}_T (\underline{u}_T), \nabla w)_T & = - (u_T, \Delta w)_T + \sum_{E_\in\ET} \omega_{TE} (u_E, \nabla w \cdot n_E)_E \quad \forall w \in \Poly{k+1}(T,\Real^2),
	\\
    \int_T r^{k+1}_T (\underline{u}_T) &= \int_T u_T.
\end{aligned}
\end{equation}
The reconstruction is exact for polynomials up to order $k+1$.
This lifts by one order the accuracy of the method with respect to a Finite-Element or Discontinuous-Galerkin scheme with the same polynomial degree.

On the other hand, a hybrid space allows the design of a scheme where element unknowns are coupled exclusively with degrees unknowns associated with the same element;
as a result, element unknowns can be expressed in terms of edge unknowns and be eliminated from the global system through static condensation (see~\cite[Appendix B]{dipietroHybridHighOrderMethod2020}).
The space for pressure variables is
\[
P^k_{h} = \left\{ q_h\in L^2(\Omega):  q_h|_T \in \Poly{k}(T)\ \forall T\in\T_h \right\}.
\]
To take into account boundary conditions and average constraint, introduce the subspaces
\[
\underline{U}^k_{h,0} \coloneq \left\{ \underline{v}_h \in\underline{U}_h^k: v_E = 0 \, \forall E \subset\partial\Omega \right\},
\quad
P^k_{0,h} \coloneq \left\{ q_h\in P^k_h: \int_\Omega q_h = 0 \right\}.
\]
\partitle{Discrete Problem}
Introduce the bilinear forms $~{a_h:\underline{U}^k_{h,0}\times\underline{U}^k_{h,0}\to\Real}$ and $~{b_h:\underline{U}^k_{h,0}\times P^{k}_h}$ such that
\[
\begin{aligned}
	a_h (\underline{u}_h, \underline{w}_h) &= \sum_{T\in\Th} \left\{ (\mu \nabla_s r^{k+1}(\underline{u}_T), \nabla_s r^{k+1}(\underline{w}_T))_T + s_T (\underline{u}_T, \underline{w}_T) \right\}\\
	b_h(\underline{w}_h, q_h) &= \sum_{T\in\Th} - \left\{ (w_T, \nabla q_T)_T - \sum_{E\in\ET} \omega_{TE} (w_E\cdot n_{TE}, q_T)_E\right\},
\end{aligned}
\]
where the $s_T:\underline{U}_T^k\times\underline{U}_T^k\to\Real$ is a positive symmetric stabilization term.
The discrete problem reads:
find $(\underline{u}_h, p_h)\in \underline{U}^k_{h,0} \times P_{h,0}^k$ such that:
\begin{equation}\label{eq:stokes.discrete}
\begin{aligned}
	a_h (\underline{v}_h, \underline{w}_h) + b_h(\underline{w}_h, p_h) &= -\sum_{E\in\Gamma} (f_\Gamma,w_E)_E & \quad &\forall \underline{w}_h \in \underline{U}^k_{h,0}, \\
    b_h(\underline{w}_h, q_h) &= 0 &\quad &\forall q_h \in \Poly{k}(\Th).
\end{aligned}
\end{equation}
\begin{remark}[Case $k=0$]
    When $k=0$, the stability of the scheme\ref{eq:stokes.discrete} is no more guaranteed.
    A solution inspired by [\cite{bottiLoworderNonconformingMethod2019}] consists in adding a symmetric penalty term that penalizes the jump of reconstructed velocity at element-element interfaces~[\cite[Section 7.6]{dipietroHybridHighOrderMethod2020}].  
    This technique is applied in the numerical tests presented in Section~\ref{sec:num-tests}.
\end{remark}
\subsection{Membrane Position Update}\label{subsec:simul.algorithm}
Let $\tau$ be a time step.
Given the total number of time iterations $N$, we introduce a sequence (times series) of polygons $\{\Gamma^{(n)}\}_{n\in\{0,..,N\}}$ such that $\Gamma^{(n)}$ approximates $\gamma(n\tau)$.
Given a discrete interface $\Gamma$ and a fixed background mesh $\Tho$ of $\Omega$, let $(\Th,\Eh)=\mathcal{M}(\Tho,\Gamma)$ be the (polygonal) mesh compatible with $\Gamma$ in the sense described in Section~\ref{subsec:hho-stokes} and illustrated in Fig.~\ref{fig:compatible.mesh} (see Appendix~\ref{app.sec:polymemb} for the details on the implementation). 
Given a compatible mesh $(\T_h, \E_h)$, let $(p_h, \underline{v}_h) = \mathcal{S}((\T_h,\E_h); f_\Gamma)$ be the corresponding solution to Problem~\eqref{eq:stokes.discrete}.
Since every edge on $\Gamma$ coincides with an edge in $\E_h$, it is possible to extract an interface velocity $v\in\Poly{k}(\Gamma,\Real^2)$ such that, for any $E$ in $\E_h$, its restriction $v_E$ is taken as the corresponding edge velocity variable of $\underline{v}_h$.
Introduce the projected velocity $\tilde{v} = \pi_{\Pone{\Gamma,\Real^2}} (v)$.
Being taken in $ \Pone{\Gamma,\Real^2}$, $\tilde{v}$ is determined by its values at the vertices.
In particulat, for $k=0$, $\tilde{v}_V = \frac12(v_{E^{(-)}}  + v_{E^{(+)}})$.
A vertex-based interface velocity allows to define a predictive advanced position $\Gamma^*$ by advancing vertices with the scheme
\[
\Gamma^{*}_V = \Gamma^{(n)}_V + \tau \tilde{v}_V \qquad \forall V\in \V_\Gamma.
\]
The membrane flow is incompressible. As a consequence, at the continuous level, the volume enclosed by the membrane $\gamma$ is expected to be conserved.
Since an explicit membrane time advancement does not guarantee this property at the discrete level, a correction step is put in place in order to enforce volume conservation.
The correction is designed in the spirit of [\cite{deliaOptimizationbasedMeshCorrection2016}] and is similar in fashion to the \textit{a posteriori} volume conservation presented in~[\cite{mendez2014unstructured}].
If $\Gamma^{n}$ is the current position of the discrete membrane, $\Gamma^*$ is the predicted configuration after advection and $V(\Gamma)$ is the volume enclosed by a polygon $\Gamma$, The corrected position configuration $\Gamma^{(n+1)}$ is found as the solution of the problem:
Find $~{\hat{\Gamma}\in \mathbb{P}^1_{(VC)}(\Gamma, \Real^2)\coloneq\left\{ \Theta \in \Pone{\Gamma^{(n)}, \Real^2}: V(\Theta) = V(\Gamma^{n})\right\}}$ such that
\begin{equation}\label{eq:min.pb}
\hat{\Gamma}  = \argmin_{\Theta\in \mathbb{P}^1_{(VC)}(\Gamma, \Real^2)} \sum_{V\in\V} \norm{}{ \Theta_V - \Gamma^{*}_V}^2.
\end{equation}
Call $\mathcal{C}(\Gamma^*,\Gamma^{(n)})$ the solution of Problem~\eqref{eq:min.pb}.
To summarize, the sequence $\{\Gamma^{(n)}\}_{n\in\{0,..,N\}}$ is determined as in Algorithm~\ref{alg:memb.update}.

\begin{algorithm}[H]\centering
\begin{algorithmic}
\Initialize{$\Gamma^{(0)}=\Gamma_0$}
\For{$ n = 0,\dots,N$}
    \State $(\Th,\Eh)^{(n)}                 = \mathcal{M}(\Tho; \Gamma^{n})$  \Comment{Generate compatible mesh}
    \State $f_\Gamma                        = f_\Gamma(\Gamma^n)$             \Comment{Compute discrete force density}
    \State $(p_h, \underline{v}_h)          = \mathcal{S}((\Th,\Eh)^{(n)}; f_{\Gamma^{(n)}})$  \Comment{Compute velocity by solving Problem~\eqref{eq:stokes.discrete}}
    \State $\Gamma^*_V                      = \Gamma^{(n)}_V + \tau \tilde{v}_V$ \Comment{Position advancement: prediction step}
    \State $\Gamma^{(n+1)}                  = \mathcal{C}(\Gamma^*, \Gamma^{n})$ \Comment{Position advancement: correction step}
\EndFor
\end{algorithmic}
\caption{polygonal membrane flow simulation}\label{alg:memb.update}
\end{algorithm}

\section{Numerical Tests}\label{sec:num-tests}
The simulation is implemented using the Python package \verb|polymemb| \footnote{\href{https://aurelio-spadotto.github.io/poly-memb/}{https://aurelio-spadotto.github.io/poly-memb/}} to perform mesh clipping and the solution of~\eqref{eq:stokes.discrete}.
The package is presented in Appendix~\ref{app.sec:polymemb}.
The minimization problem~\eqref{eq:min.pb} is solved with a Sequential Least Square Programming (SLSQP) algorithm [\cite{boggsSequentialQuadraticProgramming1995}].
The Section is organized as follows: in Section~\ref{subsec:test.convergence}, we display a set of preliminary tests with the goal of assessing numerically the convergence properties of the numerical scheme~\eqref{eq:stokes.discrete}.
In Section~\ref{subsec:test.shear} the full simulation algorithm is tested on the standard benchmark of a capsule set in a pure shear flow. 

\subsection{Preliminary convergence test}\label{subsec:test.convergence}
\begin{figure}[H]
\centering
\begin{subfigure}{\textwidth}
\includegraphics[width=\textwidth]{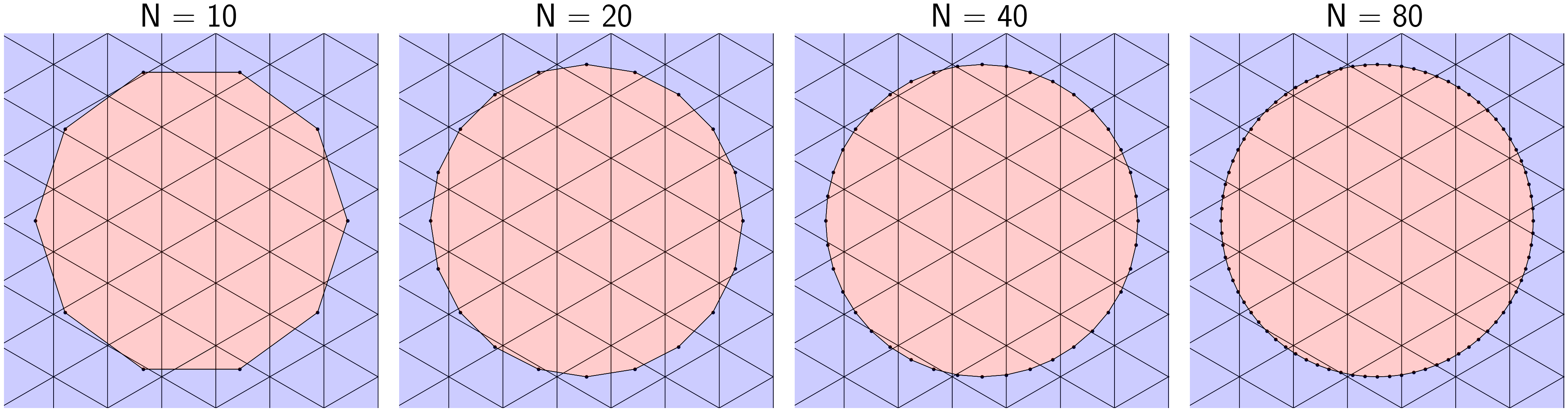}
\caption{$h_\Gamma$-convergent sequence of compatible meshes used for the convergence test. $N$ indicates the number of nodes along $\Gamma$.}\label{subfig:Gamma.sequence}
\end{subfigure}
\hfill
\begin{subfigure}{0.49\textwidth}
\includegraphics[width=\textwidth]{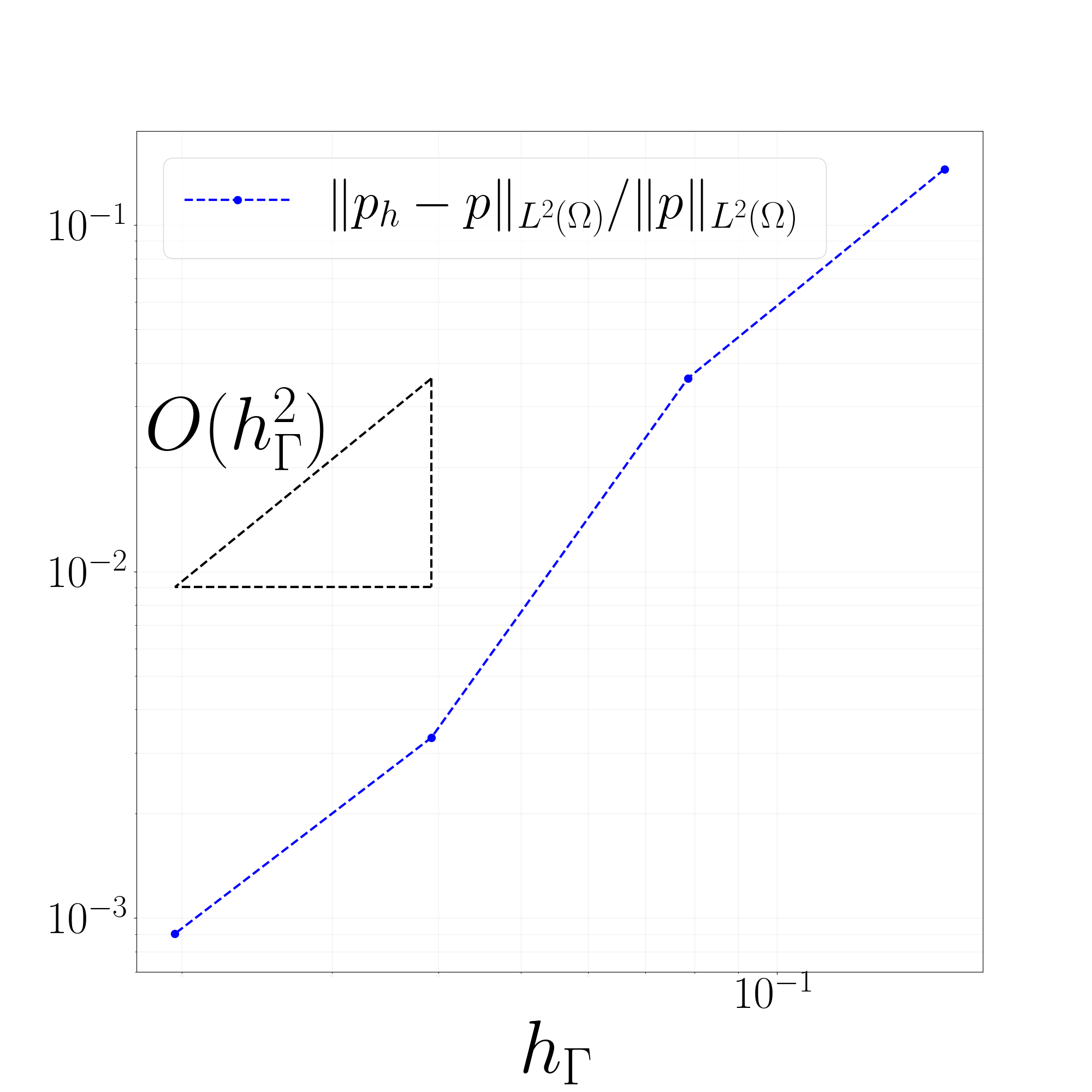}
\caption{(T1): $\lambda=1$, $\varepsilon=1e-4$}\label{subfig:conv.test.str}
\end{subfigure}
\begin{subfigure}{0.49\textwidth}
\includegraphics[width=\textwidth]{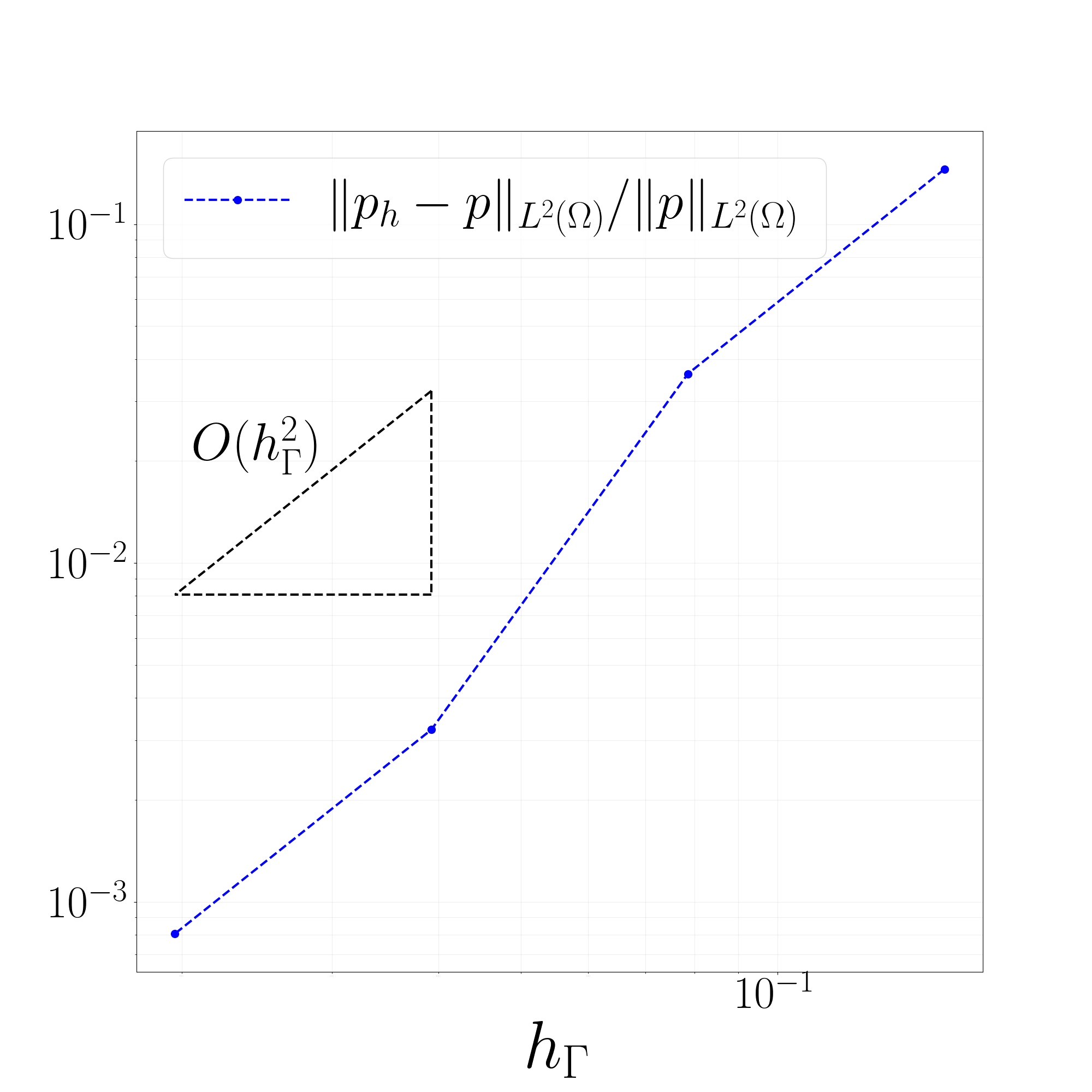}
\caption{(T2): $\lambda=1$}\label{subfig:conv.test.bnd}
\end{subfigure}
\label{fig:conv.test}
\caption{Numerical error trend for tests (T1)-(T2) from Section~\ref{subsec:test.convergence} on the sequence of progressively refined polygonal discretizations of the membrane displayed on Figure~\ref{subfig:Gamma.sequence}.}
\end{figure}

The core of the simulation is the numerical solution of Problem~\eqref{eq:stokes-pb}.
The accuracy of the approximation depends on both the interpolation properties of the discrete space $\underline{U}^k_h\times\Poly{k}(\Th)$
and the polygonal discretization of the interface $\gamma$.
For stationary configurations ($v=0$) characterized by a homogeneous $f_{\gamma}$, the solution is expressed by a region-wise constant pressure $p=(p_{\INT}, p_\EXT)$, with a constant jump across $\gamma$, reproduced exactly by Scheme~\eqref{eq:stokes.discrete}.
In such configurations, the accuracy of the solution is impacted uniquely by the geometric discretization of the interface.
With this in mind, we consider two distinct stationary configurations to assess convergence with respect to the membrane edge diameter $h_\Gamma$.
For both configurations the reference membrane $\gamma_0$ is a circle with radius $R_0=0.25$, whose center is taken as the origin of the reference system.
The two tests are designed to separate the strain and bending contribution to $f_\gamma$:
\begin{itemize}
\item[(T1)]
Pure strain response: $(\lambda_{\textrm{str}}, \lambda_{\textrm{bnd}}) = (\lambda, 0)$.
The current configuration $\gamma$ is obtained by expanding isotropically $\gamma_0$ by a factor $1+\varepsilon$.
The resulting constant force density is $f_\gamma = -\frac{\lambda\varepsilon}{R}n_\gamma$.
The corresponding solution reads
\begin{equation}
    v = 0, \qquad p_\INT - p_\EXT = \frac{\lambda\varepsilon}{R};
\end{equation}
\item[(T2)]
Pure bending response: $(\lambda_{\textrm{str}}, \lambda_{\textrm{bnd}}) = (0, \lambda)$.
We consider $\gamma = \gamma_0$.
The resulting constant force density is $f_\gamma = \frac{\lambda}{2R^3}n_\gamma$.
The corresponding solution reads
\begin{equation}
    v = 0, \qquad p_\INT - p_\EXT = -\frac{\lambda}{2R^3}.
\end{equation}
\end{itemize}
The numerical setting is common to both tests: Problem~\eqref{eq:stokes-pb} is solved in the domain $\Omega=(-1,1)^2$, discretized by a fixed background triangular mesh $\mathcal{T}_{h,0}$, with diameter $h_\Omega\approx R/2$.
The numerical pressure jump is computed selecting $k=0$ on a sequence of polygonal discretizations of $\gamma$ with edge diameter $h_\Gamma\to 0$ as illustrated in Fig.~\ref{subfig:Gamma.sequence}. 
While the interface is refined at each iteration, the background mesh is kept fixed. 
The convergence trend of the (normalized) error $\norm{L^2(\Omega)}{p-p_h}$ is reported for (T1) and (T2) in Fig~\ref{subfig:conv.test.str} and~\ref{subfig:conv.test.bnd}.
The obtained results highlight how the error introduced by the discretization of the membrane decays approximately of order 2 with respect to the edge diameter.

\subsection{Capsule in Shear Flow}\label{subsec:test.shear}

\begin{figure}[H]
\centering
\includegraphics[width=0.5\textwidth]{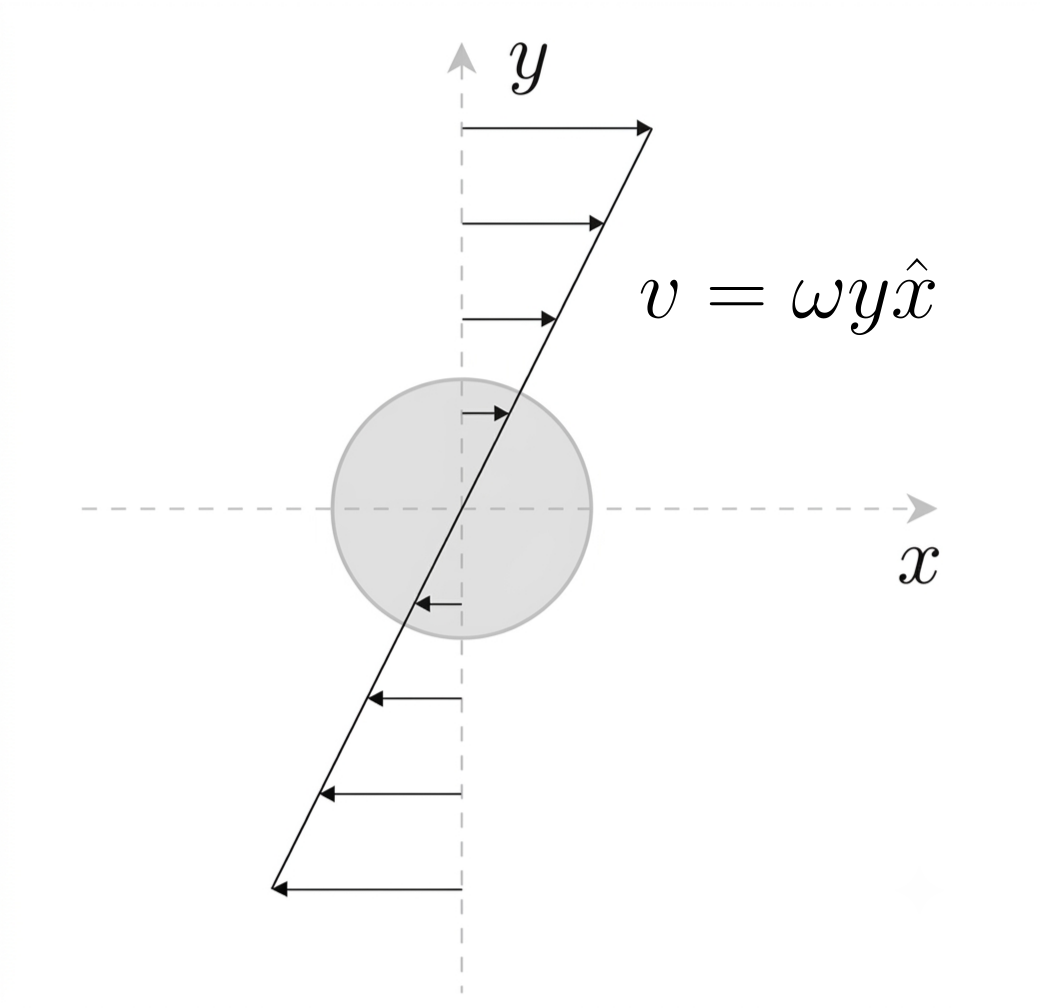}
\caption{Section~\ref{subsec:test.shear}: initial configuration of a circular capsule set in pure shear flow}\label{subfig:sketch_shear}
\label{fig:shear.overview}
\end{figure}

\begin{figure}[H]
\centering
\begin{subfigure}{0.49\textwidth}
\includegraphics[width=\textwidth]{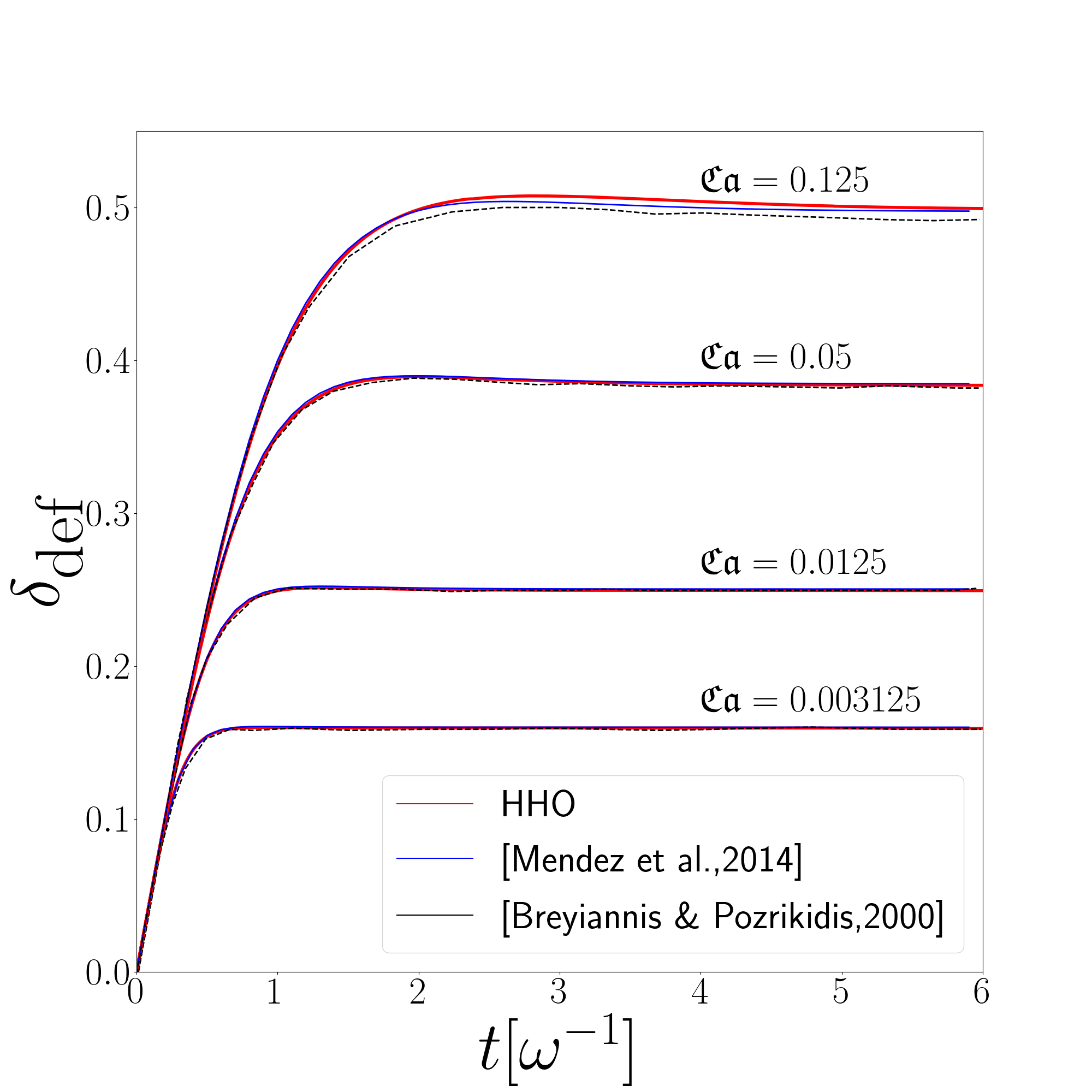}
%\caption{}\label{subfig:comparison-deformation}
\end{subfigure}
\hfill
\begin{subfigure}{0.49\textwidth}
\includegraphics[width=\textwidth]{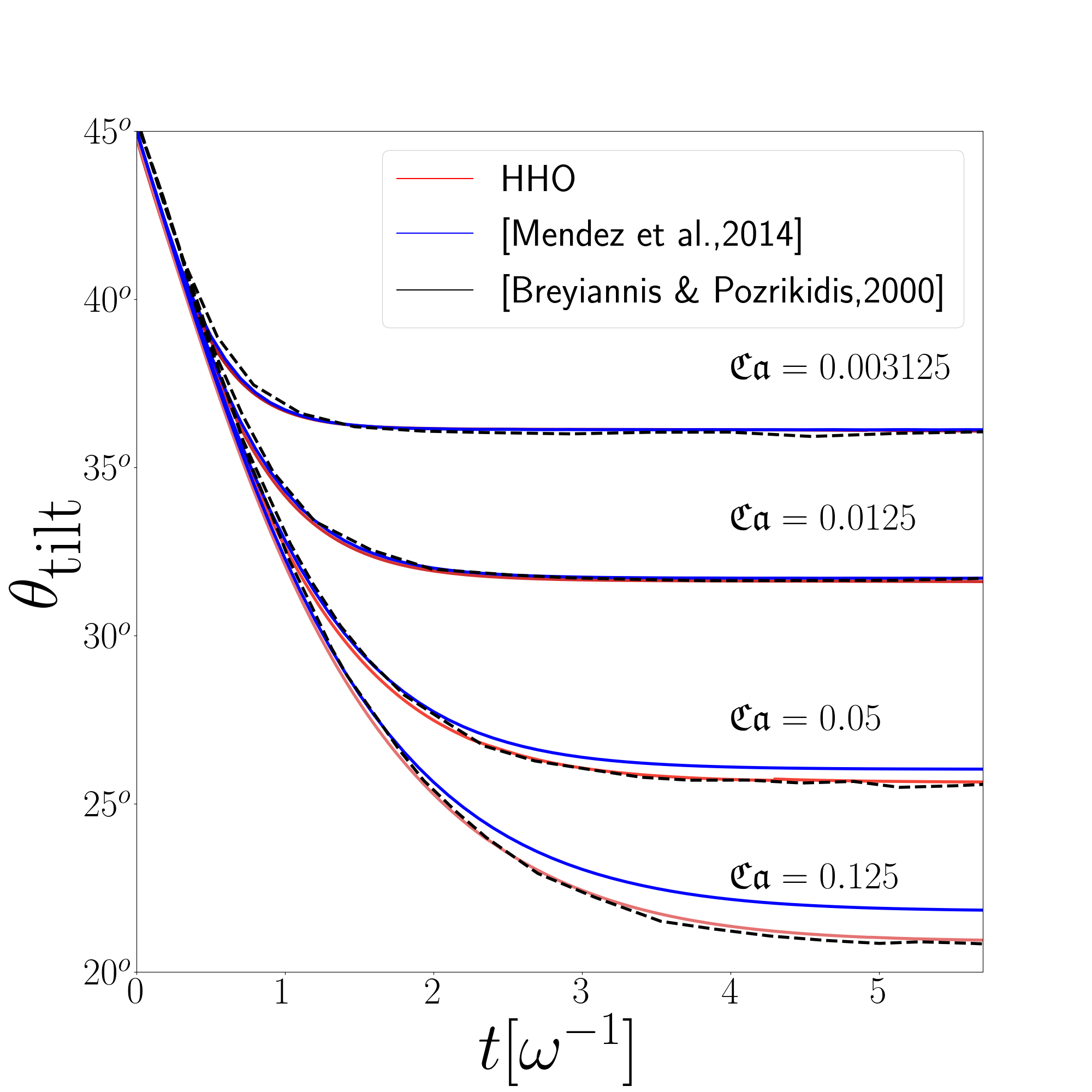}
%\caption{}\label{subfig:comparison-angle}
\end{subfigure}
\begin{subfigure}{0.49\textwidth}
\includegraphics[width=\textwidth]{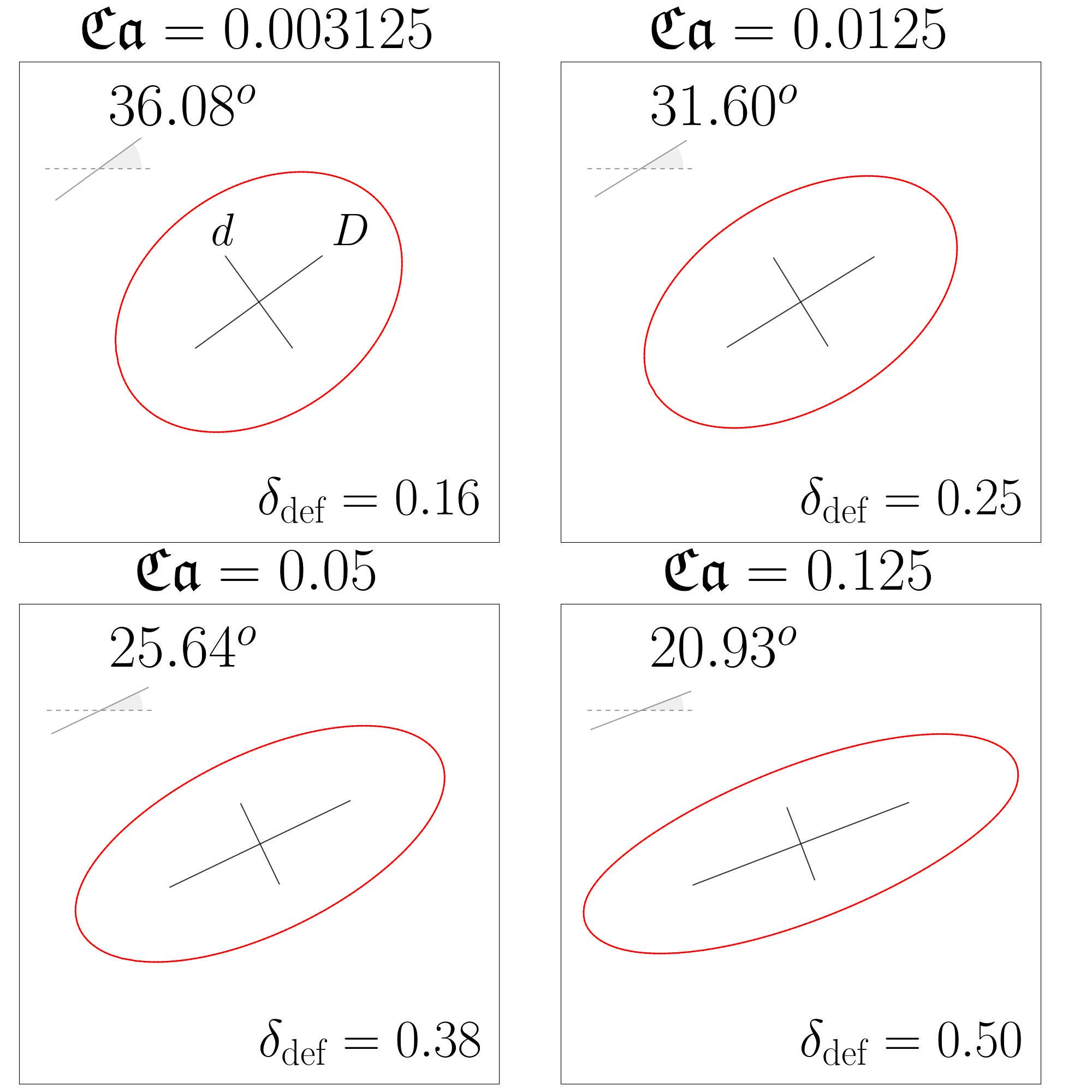}
%\caption{}\label{subfig:shear-final-state}
\end{subfigure}
\caption{
Test of Section~\ref{subsec:test.shear}: on the top, deformation $\delta_{\textrm{def}}$ and tilting angle $\theta_{\textrm{tilt}}$ obtained with the HHO-based Algorithm~\ref{alg:memb.update} compared with two simulations from literature:
an IBM-based solver applied in~[\cite{breyiannisSimpleShearFlow2000}] and Finite Volume-IBM solver adopted in~[\cite{mendez2014unstructured}].
Time is expressed with respect to the inverse shear rate~$\omega$.
At the bottom, visualization of the corresponding maximal deformation and tilting angle for the considered $\Ca$ regimes.}
\label{fig:shear.results}
\end{figure}

The goal of the test is to reproduce the well-known dynamics of an elastic capsule in shear flow.
Capsules feature a rich variety of predictable motions patterns, making it a perfect benchmark to test the validity of capsule simulation algorithms.
Under simplified assumptions (initial spherical/circular capsule, vanishing bending modulus, same viscosity inside and outside), when a capsule initially at rest is deposited in a shear flow, it is progressively deformed, initially in the principal direction of strain, that forms an angle of $\pi/4$ with the direction of the flow. 
Deformation continues until the viscous forces are balanced by the elastic response of the membrane ([\cite{barthes-bieselMotionSphericalMicrocapsule1980}],[\cite{breyiannisSimpleShearFlow2000}]).
At that point, the body of the capsule stabilizes into a tilted elliptic shape and the membrane undergoes a revolving motion around it, referred to as \textit{tank-treading}.
The maximal extent of the deformation and the tilting angle are controlled by the non-dimensional capillary number $\Ca$, measuring the ratio between the imposed viscous stresses and the restoring elastic stresses.
At high values of $\Ca$, the elastic response is weaker, capsules get elongated and tend to align more with the main flow.
\\[1mm]
The specific configuration considered for the test is represented in Fig.~\ref{fig:shear.overview}.
It is inspired by a numerical test reproduced in [\cite{breyiannisSimpleShearFlow2000}] and is based on the following specifications. 
The undeformed configuration $\gamma_0$ of the capsule is a circle of radius $R$. No resistance to bending is assumed ($\lambda_{\rm bnd}=0$).
The capsule is set in a shear flow and the fluid has uniform viscosity $\mu$.
In a Cartesian coordinate system $(\hat{x},\hat{y})$ with center set at the center of $\gamma_0$, the velocity profile far enough from the capsule is given by
\[
v = \omega y \hat{x},
\]
where $\omega$ represents the shear rate.
In this setting,
\[
\Ca \coloneq \frac{R \mu \omega}{\lambda_{\textrm{str}} }.
\]
The initial condition is displayed in Fig.~\ref{fig:shear.overview}.
A systematic parametric study is conducted to check the output of algorithm~\eqref{alg:memb.update} on a set of increasing values of $\Ca$.
The considered regimes span a range going from low-moderate deformation of the capsule ($\Ca=0.0031125$), to fairly large deformation ($\Ca=0.125$).
The $\Ca$ regime is controlled by calibrating the intensity of the shear rate $\omega$.
The time step is set as $\tau=10^{-2}\tau_{\textrm{str}}$, where $\tau_{\textrm{str}}=$ $R\mu/\lambda_{\rm str}$ represents the characteristic time of the elastic response of the capsule.
All the simulations are run up to $t=6\omega^{-1}$, therefore sufficiently beyond the exhaustion of the transient phase.
The validity of the solver is checked by tracking two characteristic quantities:
the deformation parameter
\[
\delta_{\textrm{def}}\coloneq \frac{D-d}{D+d},
\]
where $d$ and $D$ represent respectively the minor and major axis of an equivalent ellipse sharing the same inertia tensor of the capsule;
and the tilting angle ($\theta_{\textrm{tilt}}$), formed by the major symmetry axis of the capsule with the $\hat{x}$ direction.
The results of the simulation are summarized in Fig.~\ref{fig:shear.results} and show very good agreement with previously published independent results on the same configuration:
the Integral-Boundary-Method-based solver applied in~[\cite{breyiannisSimpleShearFlow2000}] (considering an unbounded domain) and a Finite Volume-Immersed-boundary-Method solver adopted in~[\cite{mendez2014unstructured}].
In particular, the simulated tilting angle is mostly coherent with [\cite{breyiannisSimpleShearFlow2000}].
This result, however, is first and foremost explained by the fact that in~[\cite{mendez2014unstructured}] flow is modelled via full Navier-Stokes equations, instead of Stokes;
the departure from~[\cite{breyiannisSimpleShearFlow2000}] may be mostly explained by nonlinear effects at higher shear rates.

\begin{figure}[H]
\centering
\includegraphics[width=\textwidth]{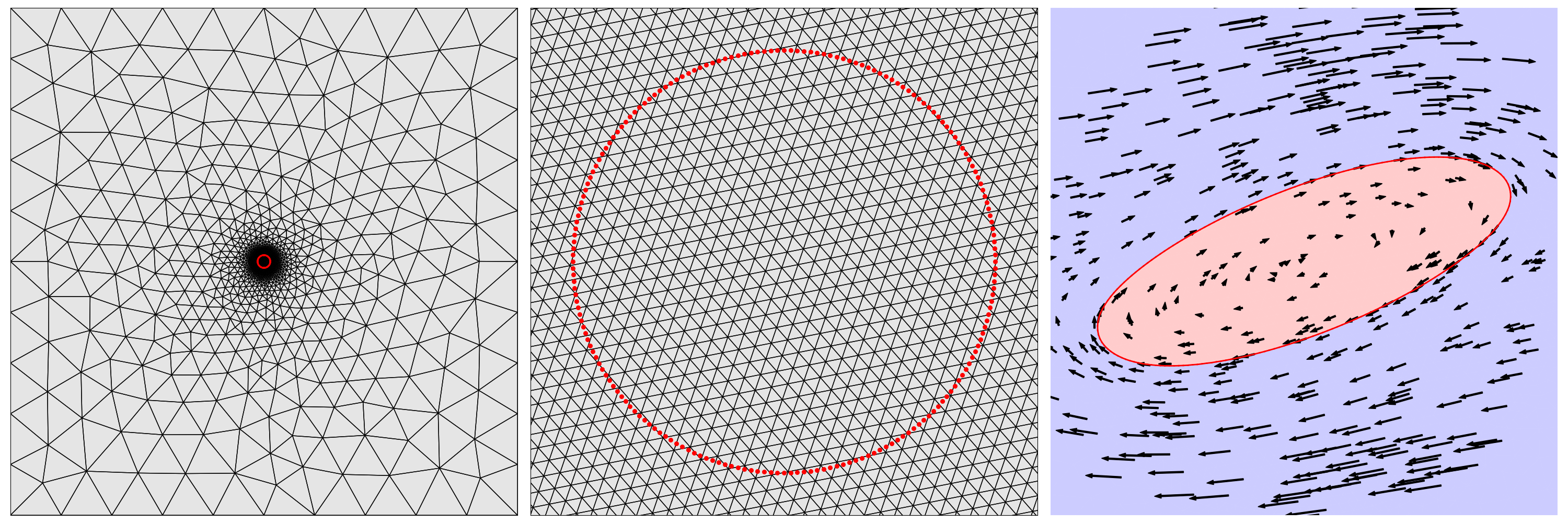}
\caption{Test of Section~\ref{subsec:test.shear}, from left to right: a) Background mesh $\Tho$ employed for the test: triangular, unstructured on the domain $\Omega=[-10, 10]^2$ with refined area close to the capsule;
b) Close-up in proximity of the interface at $t=0$ (radius $R=0.25$,  $h_\Gamma/h_\Omega \approx 1/2$);
c) Snapshot of the velocity solution for $t\gg 0$ in the regime $\Ca=0.125$.
}\label{fig:shear.mesh.setting}
\end{figure}

\begin{figure}[H]
\centering
\includegraphics[width=0.8\textwidth]{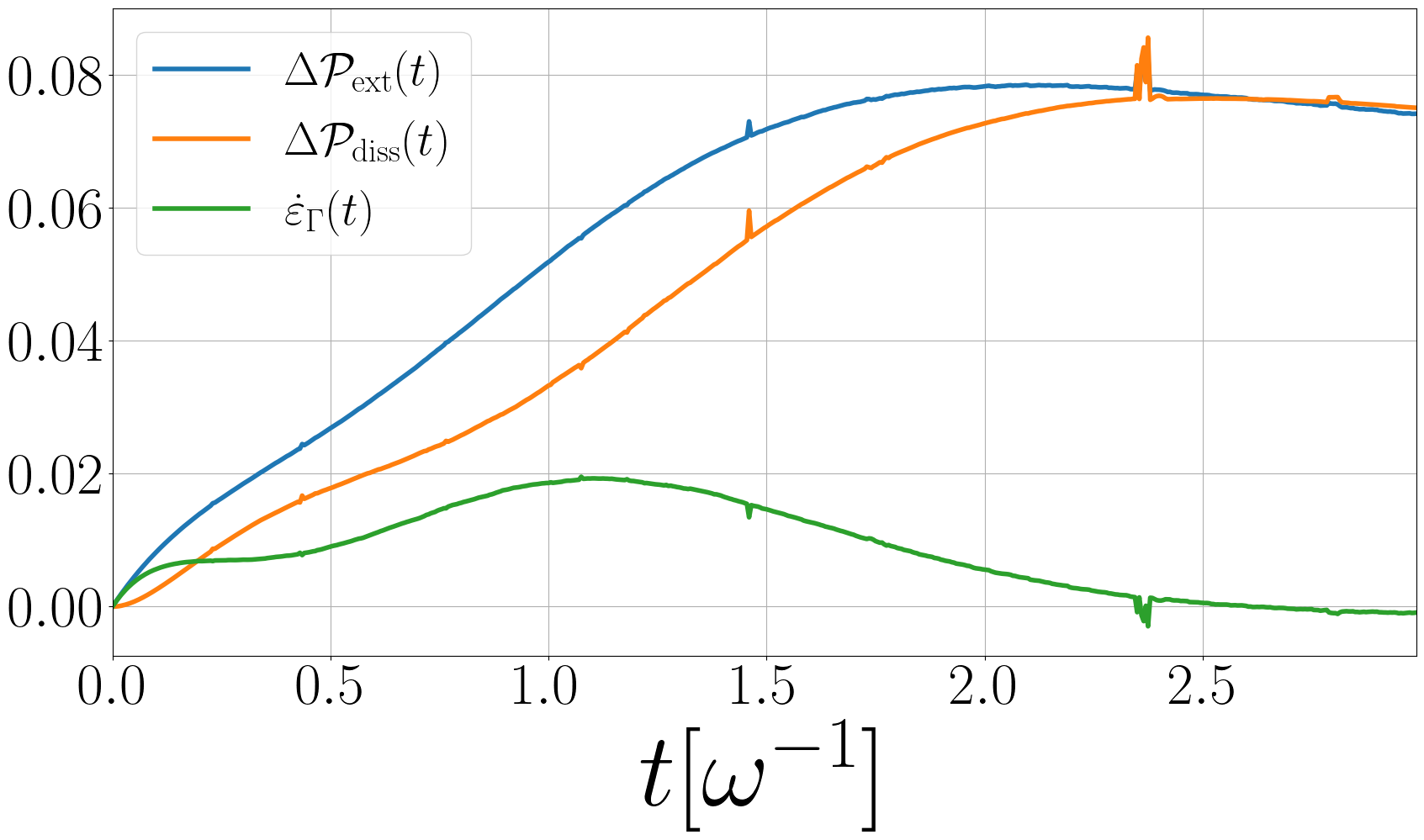}
\caption{Numerical test from Section~\ref{subsec:test.shear}, analysis of energy flows: case $\mathfrak{Ca}=0.125$; $\mathcal{P}_{\rm diss}^\infty=100$.
As pointed out in Remark~\ref{rm:sliver.cuts}, accuracy on velocity is punctually perturbed in reason of sporadical unfavorable cuts;
however, the simulation is not severely impacted at the global level.
}\label{fig:energy.balance}
\end{figure}

The HHO scheme~\eqref{eq:stokes.discrete} introduces a discretization error on velocity expressed in a norm controlling derivatives up to first order~[\cite[Lemma 8.15, Theorem 8.18]{dipietroHybridHighOrderMethod2020}].
As a consequence, it guarantees arbitrary-order accuracy in the approximation of the velocity gradient, and, therefore, of the hydrodynamical stress tensor $\sigma$.
This overcomes a known limit of the Immersed Boundary Method, which introduces a strong discretization error on the velocity gradient ([\cite{peskinImmersedBoundaryMethod2002}],[\cite{Mittal2005}]).
This property is exploited to perform a numerical verification of the energy exchanged by the fluid/capsule system.
Energetic considerations are an important tool to model the rheology of capsule suspensions;
evaluating the energy spent to induce a shearing flow to a capsule/fluid system is instrumental in estimating the effective viscosity of the suspension ([\cite{Batchelor1970, BarthesBieselRallison1981, RamanujanPozrikidis1998, GrossKrugerVarnik2014}]).  
Considering a flow with no suspended capsules, the energy injected to the control volume $\Omega$ is entirely dissipated by viscous friction at a rate
\[
\mathcal{P}_{\rm diss} = \int_\Omega \mu \nabla_s v : \nabla_s v.
\]
For a stationary flow with uniform shear rate $\omega$ this amounts to a constant dissipated power
\[
\mathcal{P}_{\rm diss}^\infty = \mu \omega^2 |\Omega|.
\]
When a capsule is present, it is deformed by the flow and stores elastic energy in the membrane at a rate: 
\[
\dot{\varepsilon}_\Gamma = \int_\Gamma v \cdot f_\Gamma.
\]
The total excess power
\[
\Delta \mathcal{P}_{\rm ext} = \Delta \mathcal{P}_{\rm diss} + \dot{\varepsilon}_\Gamma,
\qquad
\Delta \mathcal{P}_{\rm diss}= \mathcal{P}_{\rm diss} - \mathcal{P}_{\rm diss}^\infty
\]
measures the difference between a fluid/capsule and a fluid-only system in terms of total external power injected in the control volume $\Omega$.
In Fig.~\ref{fig:energy.balance}, the three quantities $\Delta\mathcal{P}_{\rm ext}$, $\Delta \mathcal{P}_{\rm diss}$ and $\dot{\varepsilon}_\Gamma$ are plotted as a function of time along the transient phase from the initial circular configuration up to the stationary tank-treading motion.
Exploiting the capacity of scheme~\eqref{eq:stokes.discrete} to compute accurate velocity gradients, this study provides, to the best of our knowledge, the first numerical verification of the dynamic energy rate balance for a single suspended capsule via direct evaluation of viscous dissipation.
The energy profile matches the expected behaviour.
The shearing flow deforms the capsule membrane, which stores elastic energy as the capsule evolves to a stable tank-treading configuration.
The presence of the capsule increases the resistance to the flow, resulting in enhanced viscous dissipation.
Initially equal to the dissipation of the capsule-free system, it increases as the capsule perturbs the flow;
after reaching a maximum during the transient phase,
viscous dissipation decreases slightly at the equilibrium, when the capsule, aligning with the direction of the flow, offers reduced resistance.
\begin{remark}[Unfavorable cuts]\label{rm:sliver.cuts}
The trajectory of the capsule is accurately reproduced without any special treatment for generated polygonal element with poor quality.
Sliver cuts are, however, the crucial challenge in the design of unfitted high-order schemes such as Cut-FEM or Cut-HHO,
which typically rely on specific techniques, such as Ghost Penalty, to enforce stability~[\cite{burmanCutFEMDiscretizingGeometry2015,burmanCutCellHybrid2019}].
Therefore, nothing prevents scheme~\eqref{eq:stokes.discrete}'s accuracy to be occasionally degraded in case of unfavorable cuts (see, for example, Fig.~\ref{fig:energy.balance}).
However, at least on the considered configuration, the noise on the approximation of velocity introduced by sliver cuts seems to have a low impact on the accuracy of the trajectory.
\end{remark}

\section*{Conclusion and Perspectives}
In this work, we have proposed a novel approach for the simulation of capsule dynamics based on interface-fitted polygonal meshes, summarized in Algorithm~\ref{alg:memb.update}.
The methodology leverages a Hybrid High Order discretization for flow equations supporting polygonal meshes.
The flexibility of the method in terms of mesh requirements allows to treat in a straightforward way the interface conditions associated to the fluid-structure interaction  while keeping the overhead linked to mesh adaptation moderate.
As a proof of concept, the performance of the polygonal interface-fitted approach has been tested successfully on a simple yet paradigmatic benchmark featuring a capsule undergoing moderate to strong deformation in shear flow.

The generality of the proposed approach goes beyond the realm of capsule mechanics, and it is possible to foresee applications to more general configurations featuring elastic shells.
An extension is the representation of interfacial forces of electric origin according to the Leaky Dielectric Model~[\cite{savilleElectrohydrodynamicsTaylorMelcherLeaky1997}].
The computation of such forces relies on the solution of an elliptic interface problem for which a strategy in the spirit of this work has been presented in~[\cite{dipietroDiscreteRhamDiscretization2025}] and implemented in \verb|polymemb|.
A limitation of the presented approach is the fact that it relies on an explicit, iteration-dependent meshing of the fluid domain.
On the one hand, this constitutes a challenge when addressing the extension to a full Navier-Stokes simulator in which
the velocity advancement scheme requires access to the fluid velocity across several time steps.
On the other hand, as illustrated in Fig.~\ref{subfig:Gamma.sequence}, an algorithm based on a polygonal solver allows to refine the interface independently from the background mesh.
This opens the way to simulations in which small-scale structures and wrinkles are fully resolved on the capsule's membrane,
while the background mesh is left relatively coarse, or locally refined in a nonconformal fashion.
Future research will explore this path to perform simulations in contexts where the membrane undergoes buckling,
such as in the case of capsules subject to compression~[\cite{barthes2011modeling}].

\section*{Acknowledgements}

Funded by the European Union (ERC Synergy, NEMESIS, project number 101115663).
Views and opinions expressed are however those of the authors only and do not necessarily reflect those of the European Union or the European Research Council Executive Agency. Neither the European Union nor the granting authority can be held responsible for them.
Simulations were performed using computer resources from the Institut de Sciences des Données de Montpellier (ISDM).
\printbibliography
\appendix
\section{Iteration-wise generation of interface-compatible meshes}\label{app.sec:polymemb}

\begin{figure}[H]
    \centering
    \includegraphics[width=\textwidth]{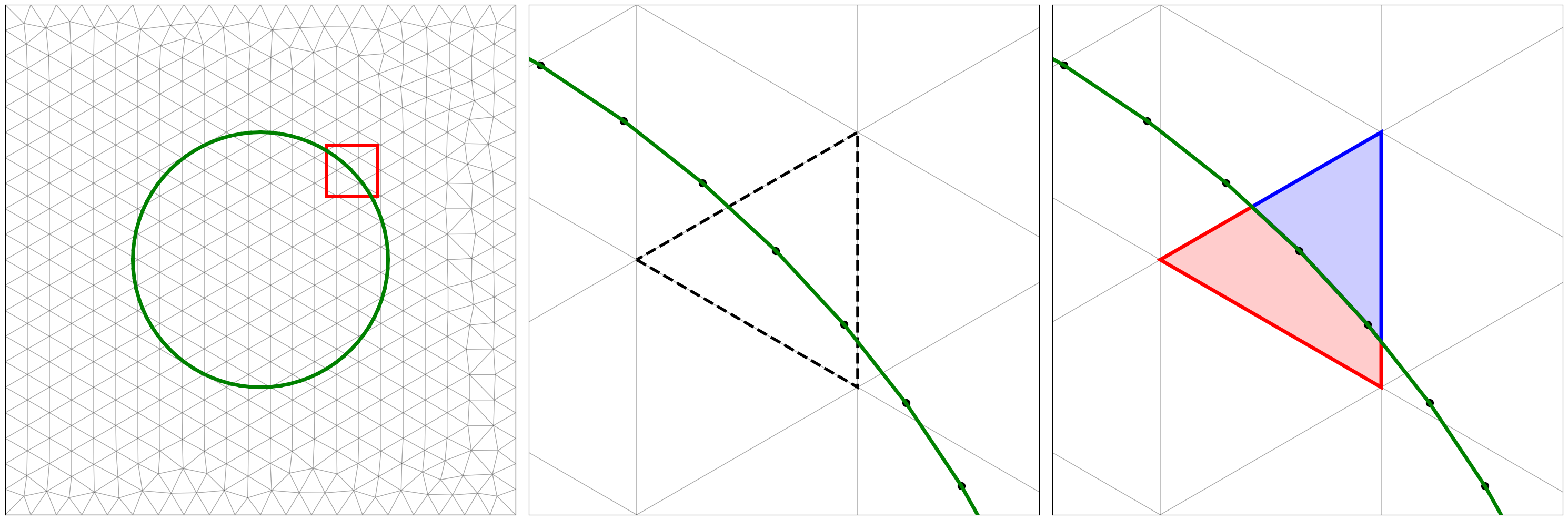}
    \caption{
    On the left, a reference simplicial mesh of a square domain and a polygonal chain (green).
    In the center, a close-up of the region framed in red, featuring a triangular element (dashed) intersected by the chain.
    On the right, the polygonal child elements (red and blue) generated by \textit{polymemb.clip} that replace the cut triangular element of the reference mesh.  
    }
    \label{fig:breakmesh-scheme}
\end{figure}

\begin{figure}[H]
    \centering
    \includegraphics[width=\textwidth]{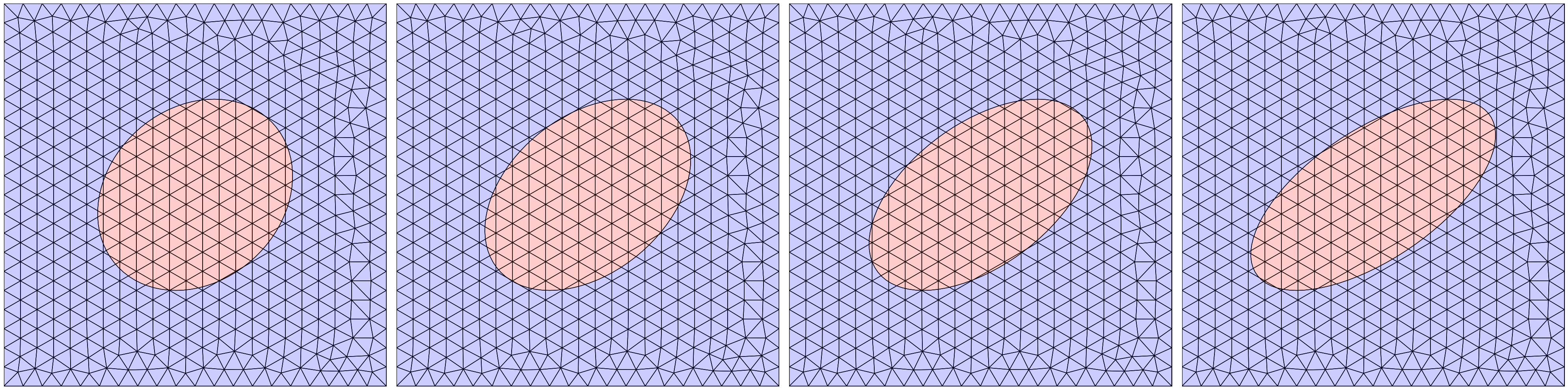}
    \caption{
    Sequence of four time snapshots and the corresponding temporary meshes generated by the clipping algorithm.
    Internal elements are marked in red, external elements in blue.
    The same background simplicial mesh is used as primer to generate each temporary mesh.
    }
    \label{fig:breakmesh-sequence}
\end{figure}

A crucial step of Algorithm~\ref{alg:memb.update} is the determination of a mesh compatible with the current configuration of the interface.
Since the interface moves across iterations, the procedure requires the generation of a new compatible mesh at each time step.
In all the numerical tests that have been presented, compatible meshes are generated by clipping a fixed background mesh along the current configuration of the interface.
The simplicity of this approach is justified by the fact that the HHO scheme introduced in Section~\ref{subsec:hho-stokes} supports polygonal elements,
which are an inevitable product of clipping.
The goal of this Section is to cover implementation aspects related to the mesh clipping procedure
at the base of the Python package \verb|polymemb|, used to produce the numerical simulations presented in this work.
\\[1mm]
The library integrates a function \verb|polymemb.clip|, which, given a generic mesh of $\Omega$ and a polygonal chain, returns a new mesh where all the elements cut by the chain are replaced by a couple of polygonal elements lying on either side (Fig.~\ref{fig:breakmesh-scheme}).
The simulation is initialized by loading a background mesh $\mathcal{T}_{h,0}$ of $\Omega$ that is kept as a reference during the entire execution of the program.
At each iteration, once the position of the interface has been updated, \verb|polymemb.clip| is executed to generate a compatible mesh (Fig.~\ref{fig:breakmesh-sequence}).
In \verb|polymemb|, members of the class \verb|mesh| are described by the list of node coordinates \verb|mesh.coords| and the list \verb|mesh.elem2node| encoding the connectivity.
The $k$-th element of \verb|mesh.elem2node| is a list of integers corresponding to the global indices of the vertices of the $k$-th element.
Element vertices are listed counter-clockwise.
\begin{lstlisting}[language=Python]
class mesh:
    def __init__(self, coords, elem2node, ...):
        """
        Args:
           coords (list[np.ndarray])  : node coordinates
           elem2node (list[list[int]]): element-to-node connectivity
        """
        self.coords = coords
        self.elem2node = elem2node
        ...
    ...
\end{lstlisting}
Polygonal chains are represented as instances of the class \verb|interface|.
The relevant member of the class is the coordinate list \verb|interface.coords|,
listing counter-clockwise the vertices.
\begin{lstlisting}[language=Python]
class interface:
    def __init__(self, coords, ...):
        """
        Args:
           coords (list[np.ndarray])  : node coordinates
        """
        self.coords = coords
        ...
    ...
\end{lstlisting}
Iteration-dependent polygonal meshes generated by \verb|polymemb.clip| are described as instances of the special class \verb|clipped_mesh|.
Clipped meshes retain a reference to a \verb|mesh| corresponding to the starting background mesh and the \verb|interface| along which clipping is performed.
For a clipped mesh, it is important to establish the connectivity between mesh edges lying along the cut and the corresponding edges on the \verb|interface|.
This is necessary to transmit the information on membrane force stored on the interface to the fluid dynamics solver, and, conversely, to broadcast flow velocity to the interface.
This information is encoded in the list \verb|clipped_mesh.intf_edges|.
\begin{lstlisting}[language=Python]
class clipped_mesh:
    def __init__(self, mesh_0, intface, ...):
        """
        Args:
           mesh_0  (polymemb.mesh)   : fixed background mesh
           intface (polymemb.intface): polygonal interface along which to clip
        """
        self.mesh_0 = mesh_0
        self.intface = intface
        [
        self.mesh,            # new polymemb.mesh instance
        self.intf_edges    # mesh-to-intface connectivity 
        ] = polymemb.clip (mesh_0, intface)   #clipping function
        ...
    ...
\end{lstlisting}
Each element of the \verb|intf_edges| is a dictionary, collecting the indices of the two elements sharing the edge lying respectively in $\Omega_{\INT}$ and $\Omega_{\EXT}$.
The dictionary also takes note of the index of the edge on the \verb|interface|, and other complementary information to reconstruct directly the mesh/interface connectivity.

\begin{figure}[H]
  \centering
  \begin{tabular}{ c  c  }

    \begin{minipage}{0.6\textwidth}
      \includegraphics[width=\linewidth]{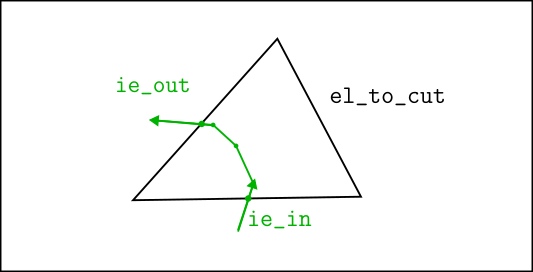}
    \end{minipage}
    \vspace{1mm}
    &
    \begin{minipage}{0.4\textwidth}
    1) The edges of the interface are looped over in order until the first that injects an element of \verb|mesh_0| (\verb|el_to_cut|) is found and flagged as \verb|ie_in|.
    The counter on the edges of \verb|intface| is increased until the edge getting out of \verb|el_to_cut| is found and flagged as \verb|ie_out|.
    The section of \verb|intface| cutting \verb|iel_to_cut| is determined.
    
    \end{minipage} 
    \\
    \begin{minipage}{0.6\textwidth}
      \includegraphics[width=\linewidth]{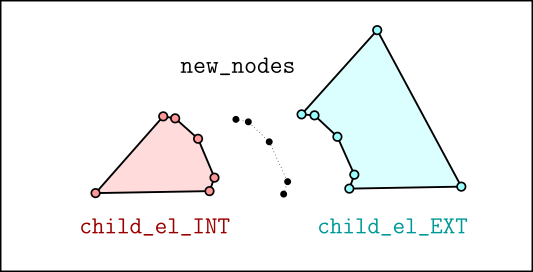}
    \end{minipage}
    \vspace{1mm}
    &
    \begin{minipage}{0.4\textwidth}
    2) The couple of new compatible elements is assembled.
    The coordinates of the intersection point of \verb|iel_in| and \verb|iel_out| with \verb|el_to_cut| and of the collection of interface nodes falling inside \verb|el_to_cut| is collected in the list \verb|new_points| and indexed.
    Two elements \verb|child_el_in| and \verb|child_el_out| are assembled as lists of node indexes and indexed. 
    \end{minipage} 
    \\
    \begin{minipage}{0.6\textwidth}
      \includegraphics[width=\linewidth]{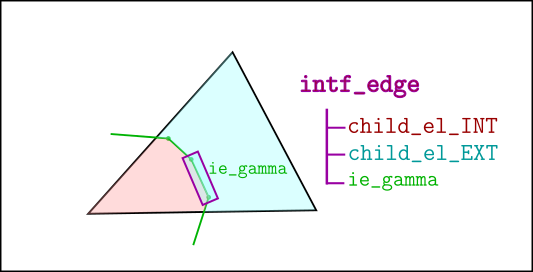}
    \end{minipage}
    \vspace{1mm}
    &
    \begin{minipage}{0.4\textwidth}
    3) For every element edge of the cut, a dictionary of type \verb|intf_edge| stores the indexes of the generated elements,
    the local index of the edge with respect to both elements and the corresponding index on \verb|intface|.
    \end{minipage} 
    \\
    \begin{minipage}{0.6\textwidth}
      \includegraphics[width=\linewidth]{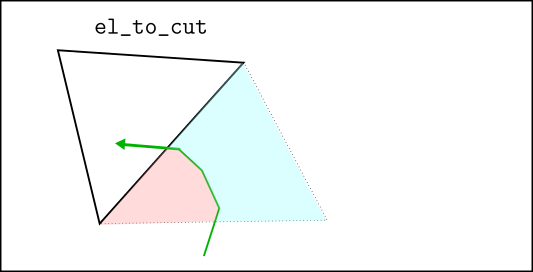}
    \end{minipage}
    \vspace{1mm}
    &
    \begin{minipage}{0.4\textwidth}
    The next element to cut is found and steps 2) and 3) are repeated on the new \verb|el_to_cut|.
    The procedure is repeated until the interface has been completely explored.
    \end{minipage} 
    \\
  \end{tabular}
  \caption{Step-by-step illustration of the implementation of the clipping process.}\label{fig:polymemb.clip}
\end{figure}

\begin{lstlisting}[language=Python]
intf_edge =  {    # dictionary
    "elem_int"     : ..., # (int): internal adjacent element
    "elem_ext"     : ..., # (int): external adjacent element
    "local_idx_int": ..., # (int): local index on internal element
    "local_idx_ext": ..., # (int): local index on external element
    "intface_ie"   : ..., # (int): index of associated edge on the interface
    ...
}    
\end{lstlisting}
The \verb|polymemb.clip| procedure is structured as an ordered loop over the edges of the polygonal chain until the first edge entering an element of the background mesh $\mathcal{T}_{h,0}$.
When such an element is detected, the algorithm proceeds by generating the child elements, flagging the shadowed cut element and recording the mesh/interface connectivity.
Once the first element has been processed, the search goes on with the next until the chain has been completely explored.
The element clipping procedure is illustrated in Fig.~\ref{fig:polymemb.clip}.

\end{document}